\documentclass[12pt]{article}
\usepackage{amsmath,amsthm,amsfonts,amssymb,MnSymbol,mathrsfs}
\usepackage[english]{babel}

\usepackage[applemac]{inputenc}
\usepackage[all]{xy}
\numberwithin{equation}{section}

  \newcommand{\const}{\rm const}
  \newcommand{\Var}{\rm Var}
  
  \newcommand{\Law}{\rm Law}

  \newcommand{\Dom}{\rm  Dom}

  \newcommand{\Sub}{\rm Sub}
  \newcommand{\StSub}{\rm StSub}
  \newcommand{\Ent}{\rm Ent}

\theoremstyle{plain}
\newtheorem{theorem}{Theorem}[section]
\theoremstyle{theorem}

\newtheorem{definition}{Definition}[section]
\newtheorem{remark}{Remark}[section]
\newtheorem{example}{Example}[section]
\newtheorem*{thma}{Theorem 4.1a}
\newtheorem*{thma2}{Theorem 4.2a}

\renewenvironment{proof}{{\bf{Proof.}}}{\hfill $\Box$ \\}

\begin{document}

  \title{Exponential tail estimates in the Law of Ordinary Logarithm (LOL) for arrays of random variables}

   \date{}
\author{\textbf{M. R. Formica, Yu. V. Kozachenko,  E. Ostrovsky and L. Sirota}}

\maketitle

\begin{center}
 Universit\`{a} degli Studi di Napoli Parthenope, via Generale Parisi 13, Palazzo Pacanowsky, 80132,
Napoli, Italy. \\

e-mail: mara.formica@uniparthenope.it \\

\vspace{4mm}

 Department of Mathematics, Kiev State University, Ukraine. \\

 e-mails: yvk@univ.kiev.ua \\
 \hspace{0.7cm} ykoz@ukr.net \\

\vspace{4mm}

Department of Mathematics and Statistics, Bar-Ilan University, \\
59200, Ramat Gan, Israel. \\

e-mail: eugostrovsky@list.ru\\

\vspace{4mm}

\ Department of Mathematics and Statistics, Bar-Ilan University, \\
59200, Ramat Gan, Israel. \\

e-mail: sirota3@bezeqint.net \\

\vspace{4mm}

\end{center}

\begin{abstract}

 \ We derive exponential bounds for tail of distribution for natural,
i.e. under ordinary  logarithm,  normalized  sums  of {\it arrays}
of random variables, not necessarily independent.\par
\end{abstract}

 \par

\vspace{3mm}

 \section{Statement of the problem. Notations. Previous
 results.}\label{intro}

\vspace{3mm}

 \ Let $(\Omega, B, {\bf P})$ be a non-trivial suitable probability
space. Let $ \xi_i,i = 1,2,\ldots$, be a {\it sequence} of centered
(i.e. with mean zero   $ \ {\bf E} \xi_i = 0 \ $ ) independent identically
distributed (i.i.d.) random variables (r.v.) having a finite non -
zero variance
$$
\sigma^2 := \Var ( \xi_i) \in (0,\infty).
$$

Denote $S_n = \sum_{i=1}^n \xi_i$, for any $n\in\mathbf{N}$. The
classical Law of Iterated Logarithm (LIL) due  by  P. Hartman and A.
Wintner \cite{Hartman-Wintner1941} tell us that
\begin{subequations}\label{limsup}
\begin{equation}\tag{\ref{limsup}}
\overline{\lim}_{n \to \infty} \frac{S_n}{\sqrt{2 n \ln \ln n}} =
\sigma
\end{equation}
with probability one (a.e.; a.s.). More general case of {\it
sequences} of independent non-identical distributed r.v. may be
found e.g. in \cite {Braverman1991},
\cite{Buldygin-Mushtary-Ostrovsky-Pushalsky}, \cite{Ostrovsky-Sirota
Jan 2008}, \cite{Ostrovsky1999}, \cite{Sgibnev 1996}, \cite{Wittmann
1985}, as well as for the martingales in \cite{Hall-Heyde1980},
\cite{Ostrovsky 1994}, \cite{Ostrovsky 2004},
 etc.\par
Analogously
\begin{equation}\label{liminf}
\underline{\lim}_{n \to \infty} \frac{S_n}{\sqrt{2 n \ln \ln n  }} =
-\sigma
\end{equation}
\end{subequations}
Let us introduce the following  {\it finite} r.v.
$$
\theta \stackrel{def}{=} \sup_{n \ge 3}  \frac{S_n}{\sqrt{2 n \ln \ln n  }},
$$
and its correspondent {\it tail function}
$$
T_{\theta}(u) = T[\theta](u)  \stackrel{def}{=} {\bf P}(\theta \ge u), \ u \ge 3.
$$
For the {\it alternating} random variables $\theta $ the  tail
function is defined as follows

$$
T_{\theta}(u) = T[\theta](u) \stackrel{def}{=} \max \left\{ \ {\bf P}(|\theta| \ge u) \ \right\}, \ u \ge 0,
$$
 the classical definition, or

$$
T^{(B)}_{\theta}(u) \stackrel{def}{=} \max \left\{ \ {\bf P}(\theta \ge u), \ {\bf P}(\theta \le -u)   \ \right\}, \ u \ge 0,
$$
the so - called  Bernstein's version, see \cite{Bernstein1964}.\par

\vspace{4mm}

 \ The exponential bound for  this tail function, e.g. of the form
\begin{equation}\label{expbboundT}
T_{\theta}(u) \le \exp \left(- C u^{m} \ \ln^r(u)  \right), \ m =
{\const} > 0, \ r = {\const} \in \mathbf{R},  \ u \ge e,
\end{equation}
 was first obtained in \cite{Kozachenko-Ostrovsky 1985}, \cite{Ostrovsky 1994}; see also \cite[chapter 2, section 2.6.]{Ostrovsky1999}.  \par

\vspace{3mm}

{\it The situation is quite different if we consider an} {\bf array}
{\it instead of the} {\bf sequence} {\it of the r.v., see e.g.}
\cite{Hu 1991}, \cite{Hu-Weber 1992}, \cite{Qi 1994}, \cite{Stout
1974}, \cite{Sung 1996}, \cite{Teicher 1981}, \cite{Wittmann 1985}.\par

\vspace{3mm}

 \ Namely, let $ \{\xi_{n,i} \}, \ i = 1,2,\ldots,n; \ n = 1,2,\ldots \
$ be an array of independent random variables with ${\bf E}
\xi_{n,i} = 0$ and  such that $0 < {\bf E} \xi^2_{n,i} < \infty.$
Define as before
\begin{equation}\label{sum array}
S_n := \sum_{i=1}^n \xi_{n,i}, \ \ s_n^2 := \sum_{i=1}^n {\bf E}
\xi^2_{n,i}, \ \  \overline{\xi}_n := \max_{i = 1,2,\ldots,n}
\xi_{n,i},
\end{equation}
and put
\begin{equation}\label{tn}
t_n :=\sqrt{2 \ s_n^2 \ \ln s^2_n },
\end{equation}
then under appropriate conditions
\begin{subequations}\label{limsup1}
\begin{equation}\tag{\ref{limsup1}}
\overline{\lim}_{n \to \infty} \frac{S_n}{t_n} = 1 \ a.e.
\end{equation}
and symmetrically
\begin{equation}\label{liminf1}
\underline{\lim}_{n \to \infty} \frac{S_n}{t_n} = - 1 \ a.e.,
\end{equation}
\end{subequations}
Law of Ordinary Logarithm (LOL).\par

 \ Evidently, if the centered r.v. $\{ \ \xi_{n,i} \ \}$ are
independent and identically distributed (i.i.d.) with
$$
\sigma^2 := {\bf E} \xi^2_{n,i} \in (0, \infty),
$$
then
$$
t_n= \sqrt{ \ 2 \ n \sigma^2 \ \ln( n \ \sigma^2 )  \ } \asymp \sqrt{ \ 2 \ n \ \ln n  \ }, \ n \to \infty.
$$
More generally, let $z = \{ \ z_n \ \}, \ n = 1,2,\ldots$ be an arbitrary
deterministic positive increasing sequence such that $\lim_{n \to \infty} z_n = \infty;$ denote

\begin{equation}\label{Q}
Q_z(u) \stackrel{def}{=} {\bf P} \left(\sup_{n} \frac{S_n}{z_n}  > u
\right).
\end{equation}

 \ For instance, the sequence $ z_n$ may coincides with $t_n: \ z_n =
\sqrt{2 \ s_n^2 \ \ln s^2_n } = t_n.$ \par

\vspace{3mm}

 {\bf  Our goal in this report is to obtain exponential decreasing  estimates for the probability $Q_z(u), \ u \ge
 3$, as $u \to \infty, \ $ as well as for near tail probabilities, for certain suitable norming sequences
  $ \ \{z_n \} \ $, possibly, on different more fine form.} \par

\vspace{3mm}

  \ Analogous estimates in the classical LIL for real or Banach spaces
valued r.v., as well as for martingales, was obtained in
\cite{Ostrovsky 1994}, \cite{Ostrovsky 2004}, \cite{Sgibnev 1996};
see also \cite[chapter 2, section 2.6]{Ostrovsky1999}. \par

\vspace{3mm}

 \  Another statement of problem is represented in the recent article
\cite{Kievinaite - Saulys 2018}, where is described also an interest application
in an insurance. \par

\vspace{3mm}

 \section{Grand Lebesgue Spaces of random variables.}

\vspace{3mm}

{\it A classical approach.} \par

 \vspace{3mm}

We present here some known facts from the theory of one-dimensional
random variables with exponential decreasing tails of distributions
and the connections with the so-called Grand Lebesgue Spaces (GLS)
and the Orlicz exponential Spaces, see \cite[chapters
1,2]{Buldygin-Mushtary-Ostrovsky-Pushalsky}, \cite{Kolmogoroff 1929}
- \cite{Kozachenko-Ostrovsky-Sirota Oct2017}, \cite[chapter 2,
section 2.6]{Ostrovsky1999}, \cite{Ostrovsky 2004},
\cite{Ostrovsky-Sirota Oct2015}, \cite{Samko-Umarkhadzhiev},
\cite{Samko-Umarkhadzhiev-addendum}. The Grand Lebesue spaces and
several generalizations and variants of them have been widely
investigated, see e.g. \cite{Iwaniec-Sbordone 1992},
\cite{Fiorenza2000}, \cite{liflyandostrovskysirotaturkish2010},
\cite{Ostrovsky-Sirota-fundamental function GLS 2015},
\cite{caponeformicagiovanonlanal2013}, \cite{anatriellofiojmaa2015},
\cite{formicagiovamjom2015}. These spaces are of great interest for
their applications not only in statistics, in theory of random
fields, Monte-Carlo methods but also in the theory of Partial
Differential Equations (PDEs) (see e.g.
\cite{Fiorenza-Formica-Gogatishvili-DEA2018} and references therein,
\cite{fioformicarakodie2017}), in interpolation theory (see e.g.
\cite{fiokarazanalanwen2004}, \cite{fioforgogakoparakoNAtoappear}),
in topics concerning the boundedness of operators (see e.g. \cite
{Ostrovsky-Sirota-boundedeness operator bilateral GLS Oct2011},
\cite{fioguptajainstudiamath2008}, \cite{anatrielloformicaricmat2016}). \par

 \ Let $\lambda_0 \in (0,\infty]$ and let $ \phi = \phi(\lambda)$ be an
even strong convex function in $(-\lambda_0, \lambda_0)$ which takes
positive values, twice continuously differentiable; briefly $\phi =
\phi(\lambda)$ is a Young-Orlicz function, such that
\begin{equation}\label{Young-Orlicz function}
\phi(0) = 0, \ \ \phi'(0) = 0, \ \  \phi^{''}(0) \in (0,\infty).
\end{equation}
We denote the set of all these Young-Orlicz function as  $\Phi: \
\Phi = \{ \phi(\cdot)  \}. $

\begin{definition}
Let $\phi\in \Phi$. We say that the centered random variable $\xi$
belongs to the space $B(\phi)$  if there exists a constant $\tau
\geq 0$ such that
\begin{equation}\label{spaceB}
\forall \lambda \in (-\lambda_0, \lambda_0) \ \Rightarrow {\bf E}
\exp(\pm \lambda \ \xi) \le \exp(\phi(\lambda \ \tau)).
\end{equation}
 The minimal non-negative value $\tau$ satisfying \eqref{spaceB} for
any $\lambda \in (-\lambda_0, \ \lambda_0)$ is named $B(\phi)$-norm
of the variable $\xi$  and we write
\begin{equation}\label{Bnorm}
||\xi||_{B(\phi)} \stackrel{def}{=}\inf \{\tau \ge 0 \ : \ \forall
\lambda \in (-\lambda_0, \lambda_0) \ \Rightarrow {\bf E} \exp(\pm
\lambda \ \xi) \le \exp(\phi(\lambda \ \tau)) \} .
\end{equation}
\end{definition}
For instance if $\phi(\lambda)=\phi_2(\lambda) := 0.5 \ \lambda^2, \
\lambda \in \mathbf{R}$, the r.v. $\xi$ is \emph{subgaussian} and in
this case we denote the space $B(\phi_2)$ with $\Sub$. Namely we
write $\xi \in \Sub$ and
$$
||\xi||_{\Sub} \stackrel{def}{=} ||\xi||_{B(\phi_2)}.
$$
 \  It is known, see  \cite{Kozachenko-Ostrovsky 1985}), \cite{Buldygin-Mushtary-Ostrovsky-Pushalsky} that if the r.v. $\xi_i$ are
independent and subgaussian, then

$$
||\sum_{i=1}^n \xi_i||_{\Sub} \le \sqrt{\sum_{i=1}^n ||\xi_i||^2_{\Sub}}.
$$

 \ At the same inequality holds true in the more general case in the $ \ B(\phi) \ $ norm, when the function $ \ \lambda \to \phi(\sqrt{\lambda})  \ $ is convex,
 see \cite{Kozachenko-Ostrovsky 1985}. \par
 \ As a slight corollary: in this case and if in addition the r.v. - s $ \ \{\xi_i \} \ $ are i., i.d., then

$$
 \sup_{n = 1,2,\ldots}|| n^{-1/2} \sum_{i=1}^n \xi_i||B(\phi) =  ||\xi_1||B(\phi).
$$

\vspace{4mm}

\begin{definition}
 \ The centered r.v. $\xi$ with finite non-zero variance $\sigma^2 :=
\Var (\xi) \in (0,\infty)$ is said to be strictly subgaussian, and
write $ \ \xi \in \StSub$, iff
$$
{\bf E} \exp(\pm \ \lambda \ \xi) \le \exp(0.5 \ \sigma^2 \
\lambda^2), \ \lambda \in \mathbf{R}.
$$
\end{definition}

 \ For instance, every centered non-zero Gaussian r.v. belongs to the
space $\StSub.$ The Rademacher's r.v. $\xi$, that is such that
${\bf P}(\xi = 1) = {\bf P}(\xi = - 1) = 1/2$, is also strictly
subgaussian. Many other strictly subgaussian r.v. are represented in
\cite{Buldygin-Mushtary-Ostrovsky-Pushalsky},
\cite{Kozachenko-Ostrovsky-Sirota Jan2017}, \cite[chapter
1]{Ostrovsky1999}. \par

\vspace{3mm}

It is proved in particular that $B(\phi), \ \phi  \in \Phi$, equipped with the norm
\eqref{Bnorm} and under the ordinary algebraic operations, are
Banach rearrangement invariant  functional spaces, which are
equivalent the so-called Grand Lebesgue spaces as well as to Orlicz
exponential spaces. These spaces are very convenient for the
investigation of the r.v. having an exponential decreasing tail of
distribution; for instance, for investigation of the limit theorem,
the exponential bounds of distribution for sums of random variables,
non-asymptotical properties, problem of continuous and weak
compactness of random fields, study of Central Limit Theorem in the
Banach space, etc. \par

\vspace{3mm}

 \ Let $ \ g: \mathbf{R} \to \mathbf{R} \ $ be numerical valued measurable function, which can perhaps take the infinite value.
 Denote by $ \Dom[g] \ $ the domain of its finiteness:

\begin{equation}\label{Domain}
{\Dom[g]} := \{y\, : \, \ g(y) \in (-\infty, \ + \infty) \ \}.
\end{equation}

 \   Recall the definition $ \ g^*(u) \ $  of  the Young-Fenchel or
Legendre transform for the function $ \ g: \mathbf{R} \to \mathbf{R}
\ $:

\begin{equation}\label{ Definition of the Young-Fenchel transform}
g^*(u) \stackrel{def}{=} \sup_{y \in \Dom[g]} (y u - g(y)),
\end{equation}
 but we will use further the value $ \ u \ $  to be only non - negative.\par

 \ In particular, we  denote by $\nu(\cdot)$ the Young-Fenchel or
Legendre transform for the function $\phi\in \Phi$:
\begin{equation}\label{Young-Fenchel transform}
\nu(x) = \nu[\phi](x)  \stackrel{def}{=} \sup_{\lambda: |\lambda|
\le \lambda_0} (\lambda x - \phi(\lambda)) = \phi^*(x).
\end{equation}

 \ It is important to note that if the non-zero r.v. $\xi$ belongs to
the space $B(\phi)$ then

\begin{equation}\label{conditionB}
{\bf P}(\xi > x) \le \exp \left(- \nu(x/||\xi||_{B(\phi)}\right).
\end{equation}
 \ The inverse conclusion is also true up to a multiplicative constant
under  suitable conditions.\par

\vspace{4mm}

 \ Furthermore, assume that the {\it centered} r.v. $\xi$  has in some
non-trivial neighborhood of the origin finite \emph{moment generating
function} and define

\begin{equation}\label{momentfunction}
\phi_{\xi}(\lambda) \stackrel{def}{=} \max_{\alpha = \pm 1} \ \ln {\bf E} \exp (\ \alpha \lambda \ \xi \ ) < \infty, \ \lambda \in ( - \ \lambda_0, \lambda_0)
\end{equation}
for some $\lambda_0 \in (0, \ \infty]$. Obviously, the last
condition \eqref{conditionB} is quite equivalent to the well known
Cramer's one. \par

 \ We agree that $\phi_{\xi}(\lambda) := \infty$ for all the values $\lambda$ for which
\begin{equation}\label{mean infinity}
 {\bf E} \exp ( \ |\lambda| \ \xi) = \infty.
\end{equation}
The function $\phi_{\xi}(\lambda)$ introduced in
\eqref{momentfunction} is named {\it natural} function for the r.v.
$ \xi$; herewith  $\xi \in B(\phi_{\xi}) $ and moreover we assume

$$
||\xi||_{B(\phi_{\xi})} = 1.
$$

\vspace{4mm}

{\it Grand Lebesgue Spaces (GLS) approach.} \\

 \vspace{3mm}

 \ Let $ (\Omega, B, {\bf P})$ be again certain non-trivial suitable probability
space. Let $ \ b = \const > 1;  \ $  the case   $ \ b = + \infty \ $  is also not excluded. \par

\  Let also $ \ p \in[1, b) \ $ or $ \ p \in [1,b]; \ $ evidently, the last case take place iff the value $ \ b \ $ is finite
(and greatest than 1). Let
$\psi_{(b)}= \psi=\psi(p)$ be a continuous function  defined in the domain $[1, b)$ such that $\inf \psi(p) > 0$. \par

 \ We can and will suppose without loss of generality  $b = \sup \{p, \psi(p) < \infty\}$, so that
$ \Dom[\psi] = [1, b) $  or $ \Dom [\psi] = [1, b], $ of course iff $ \ b < \infty. \ $ \par
 \ When $ \ b < \infty, \ $ we define formally $ \ \psi(p) = + \infty \ $ for the values $ \ p > b. \ $ \par

 \ The set of all such functions will be denoted by $\Psi_{(b)} = \{ \psi(\cdot)\}$ or $\Psi := \Psi_{\infty}$ if $b=\infty$. \par
 \ Denote also

$$
U\Psi \stackrel{def}{=} \cup_{b \in (1,\infty)} \Psi_{(b)} \ \cup \Psi.
$$

 \vspace{3mm}

 \  Let $ \psi(\cdot) $ be some function from the set $ \ U\Psi. \ $  We define the following important auxiliary  function

\begin{equation}\label{def h}
h(p) = h_{\psi}(p) \stackrel{def}{=} p \ \ln \psi(p) , \   p \in \Dom[\psi].
\end{equation}

\vspace{4mm}

\begin{definition}
 \ The Grand Lebesgue Space (GLS) $G\psi = G\psi_{(b)}$ consists of all
the numerical valued random variables (measurable functions) $ \{ \ \zeta \ \} $
defined on our probability (measurable) space and having a finite norm
\begin{equation}\label{GLSnorm}
||\zeta|| = ||\zeta||_{G\psi} \stackrel{def}{=} \sup_{p \in \Dom[\psi] }
\left\{ \frac{\ ||\zeta \ ||_p}{\psi(p)} \right\},
\end{equation}
where $\|\zeta\|_{p}$ denotes the classical Lebesgue-Riesz
$L_p$-norm
$$\|\zeta\|_p=\|\zeta\|_{L_p(\Omega)}=\left( {\bf E}|\zeta|^p\right)^{\frac{1}{p}}=\left(\int_\Omega |\zeta(\omega)|^p\  {\bf P}(d\omega)\right)^{\frac{1}{p}}, \ \ p\geq 1.$$
\end{definition}
The function $\psi =\psi(p)$ is named ordinary {\it generating
function} for the Grand Lebesgue Space $G\psi$. \par

 \ Let $\xi$ be a random variable such that there exists $ p=\const > 1 $
so that $\ || \xi ||_p < \infty$. The {\it natural} $G\Psi$ function
$\psi_{\xi} = \psi_{\xi}(p)$ for the r.v. $\xi$ is defined by the relation
$$
\psi_{\xi}(p) \stackrel{def}{=} ||\xi||_p,
$$
with correspondent domain of definition $ \ \Dom \left[ \ \psi_{\xi} \ \right], \ $ bounded or not. \par

 \ The function $ \ \psi = \psi(p), \ $ finite at last for some value $ \ p \ $ greatest than one is said to be {\it natural}, iff
there exists a r.v. $ \ \xi = \xi(\omega) \ $ for which

$$
\psi(p) = ||\xi||_p.
$$
 \ The complete description of such functions may be found in \cite{Ostrovsky1999}, chapter 1, sections 1.1., 1.8. \par

  \ These GLS spaces are rearrangement-invariant Banach functional spaces in
the classical sense and were investigated in particular in  many
works, see e.g. \cite[chapter
1]{Buldygin-Mushtary-Ostrovsky-Pushalsky},
\cite{caponeformicagiovanonlanal2013}-\cite{formicagiovamjom2015},
\cite{Iwaniec-Sbordone 1992}, \cite{Kolmogoroff
1929}-\cite{Kozachenko-Ostrovsky-Sirota Oct2017}, \cite[chapters 1,
2]{Ostrovsky1999}, \cite{Ostrovsky-Sirota Oct2015},
\cite{Samko-Umarkhadzhiev}, \cite{Samko-Umarkhadzhiev-addendum},
\cite{anatrielloformicaricmat2016}, etc.

\begin{example}

\ {\rm Let $ \ \Omega = \{\omega \} = [0,1] \ $ equipped with ordinary Lebesgue measure $ \ {\bf P.}\ $
Introduce the r.v. $ \ \xi = \xi_{a,b}(\omega), b = \const \in (1,\infty), \ a = \const \in R \ $ as follows

\begin{equation}\label{Dom - b }
\xi = \omega^{-1/b} \ |\ln \omega|^a \ I_{(0,1/e)}(\omega),
\end{equation}
where $ \ I_A(\omega) \ $ denotes the ordinary indicator function of the (measurable) set $ \ A. \ $ \par
 \ The natural function $ \ \psi_{\xi} = \psi_{\xi}(p)  \ $ has the following form

\begin{equation}\label{psib1}
\psi_{\xi}(p) < \infty, \ p \in [1,b); \ \psi_{\xi}(p) = \infty, \ p > b;
\end{equation}

\begin{equation}\label{psib2}
\psi_{\xi}(b) < \infty  \ \Longleftrightarrow ab < - 1.
\end{equation}

 \ So, the domain $ \ \Dom[\psi] \ $  can be either closed as well as semi-open. \par
 }
\end{example}

\vspace{4mm}

\begin{example}
{\rm Define $\psi= \psi_{(b)}(p)=1, \ p \in [1,b], \ 1 < b < \infty$.

 \ One can define formally $\psi_{(b)}(p) = + \infty, \ p > b$. It is
easy to verify by virtue of Lyapunov's inequality
that the $G\psi_{(b)}$ norm of any r.v. $\xi$ is
quite equal to the classical Lebesgue-Riesz $L_b$-norm
\begin{equation}\label{GLS-Lebesgue}
||\xi||_{G\psi_{(b)}} = \|\xi\|_{L_b(\Omega)}.
\end{equation}
}
\end{example}

\vspace{4mm}

\begin{example}
{\rm Define $\psi= \psi_{(b)}(p)=(b-p)^{-\frac{1}{p}}, \ p \in
[1,b), \ 1 < b < \infty$.

\noindent  Let $\varepsilon \in (0,b-1)$ and replace $p$ with
$p-\varepsilon$ and $b$ with $p$. So we have
$\psi=\psi_{(p)}(\varepsilon)=\varepsilon^{-\frac{1}{p-\varepsilon}}$,
\ $\varepsilon\in (0,p-1)$, $p>1$, and the $G\psi_{(b)}$ norm of any
r.v. $\xi$ takes the well known form
\begin{equation}\label{GLS-Lebesgue}
||\xi||_{G\psi_{(p)}} =
\sup_{0<\varepsilon<p-1}\varepsilon^{\frac{1}{p-\varepsilon}}\|\xi\|_{L_{p-\varepsilon}(\Omega)}.
\end{equation}
}

\end{example}

\vspace{4mm}

\noindent Now we refer here some facts about these spaces used in the sequel. \par

\vspace{4mm}

 \ It is known (see \cite{Kozachenko-Ostrovsky 1985}, \cite{Kozachenko-Ostrovsky-Sirota Oct2017}  )
   that if  $\xi \ne 0$ and $\xi \in G\psi_{(b)}, $ including the case $ \ b = \infty, \ $ then
\begin{equation}\label{estimateT}
T_{\xi} ( y)  = {\bf P}(|\xi| > y)
\le \exp \left( \ - h_{\psi}^* (\ln ( y/||\xi||) )
\ \right), \ y \ge e\cdot ||\xi||.
\end{equation}

 \ Namely, let $||\xi||=||\xi||_{G\psi_{(b)} } = 1$. By means of Tchebychev - Markov inequality
$$
 T_{\xi} ( y)={\bf P}(|\xi| > y) \le  \frac{\psi^p(p)}{y^p} =\exp \left(  - p \ln y + p \ln \psi(p)    \right),
$$
and consequently
\begin{equation}\label{conseguence Tchebychev-Markov}
\begin{split}
T_{\xi} ( y) & \le   \inf_{ p \in \Dom[\psi] } \exp \left(- p \ln y +
p \ln \psi(p) \right) = \\
& \inf_{ p \in \Dom[h] } \exp \left( - p \ \ln y + h_{\psi}(p)  \right) = \exp \left( - h^*_{\psi} (\ln y ) \right), \ \ \ \  y \ge e,
\end{split}
\end{equation}

as long as $ \ \Dom[\psi] = \Dom[h]. \ $ \par

 \ More generally, if $||\xi||=||\xi||_{G\psi_{(b)}} \neq 1$, we
can consider the {\it normalized} r.v. $\displaystyle \xi_{(n)}=\frac{|\xi|}{||\xi||}$
so that $||\xi_{(n)}||=1$. Then
$$
T_{\xi} ( y)={\bf P}(|\xi| > y)={\bf P}\left(\frac{|\xi|}{\|\xi\|}>\frac{y}{\|\xi\|}\right)={\bf P}\left(|\xi_{(n)}|>\frac{y}{\|\xi\|}\right)
$$
and, for the previous conclusion, we obtain \eqref{estimateT}.

 \ Conversely, the last inequality may be reversed  in the following
version: if the r.v. $\xi $ {\it satisfies the Cramer's condition} and
$$
{\bf P}(|\xi| > y) \le \exp \left(-h_{\psi}^* (\ln (y/K) \right),
\  \ y \ge e \cdot K, \ \ K = \const  >0
$$
for some generating function $ \ \psi(\cdot) \ \in U\Psi, \ $
and if the function $ h_{\psi}(p), \ 1 \le p < \infty  \ $  is
positive, continuous, convex and such that
$$
\lim_{p \to \infty}  \psi(p)/p = 0,
$$
then $ \xi \in G\psi$. Furthermore there exist $C_2(\psi) >
C_1(\psi)>0$ such that
$$ C_1(\psi) K\leq ||\xi||_{G\psi}\leq C_2(\psi)K. $$

 \ More generally, let $ \ V = V(x)\ x \ge 0 \ $ be some tail function (or some its majorant), and let $ \ \xi \ $ be any r.v. such that
$ \ T_{\xi}(x) \le V(x), \ x \ge 0. \ $   As long as

$$
||\xi||^p_p = p \int_0^{\infty} x^{p-1} \ T[\xi](x) \ dx,
$$
we conclude

$$
||\xi||_p \le \left[ \  p \ \int_0^{\infty} x^{p-1} \ V(x) \ dx \ \right]^{1/p}.
$$

 \ For instance, let

$$
T_{\xi}(x) \le  T^{(\beta, \gamma, L)}(x),  \beta = \const \in
(1,\infty), \ \gamma = \const \in \mathbf{R},
$$
where by definition

$$
T^{(\beta, \gamma, L)}(x) \stackrel{def}{=} x^{-\beta} \ (\ln x)^{\gamma} L(\ln x), \ x \ge e,
$$
and

$$
\psi^{(\beta,\gamma,L)}(p) := (\beta - p)^{-(\gamma + 1)/\beta} \ L^{1/\beta} \left( \  \frac{1}{\beta - p}  \  \right), \ 1 \le p < \beta,
$$
where in turn   $ \ L = L(x), \ x \ge 1 \ $ is some positive
continuous {\it slowly varying} function as $ x \to \infty \ $; the
set of all such functions will be denoted by $ \ SV; \ SV = \{
L(\cdot) \}.  \ $ We have

$$
T_{\xi}(x) \le  T^{(\beta, \gamma, L)}(x) \ \Rightarrow ||\xi||G\psi^{(\beta,\gamma,L)} = C_1(\beta,\gamma,L) < \infty.
$$

 \ Inversely, if $ \ ||\xi||G\psi^{(\beta,\gamma,L)} = C_2 < \infty, \ $ then

$$
T{\xi}(x) \le  C_3(\beta,\gamma,L) \cdot T^{(\beta, \gamma + 1, L)}(x);
$$
and both these estimates are non - improvable, see \cite{Kozachenko-Ostrovsky-Sirota Jan2017}. \par

\vspace{5mm}

 \ Let us introduce the following {\it exponential} Young - Orlicz function
$$
N_{\psi}(u) = \exp \left(h_{\psi}^* (\ln |u|) \right),  \ |u| \ge 1;
\ \ N_{\psi}(u) = C u^2, \ |u| < 1.
$$
and we denote the correspondent Orlicz norm by $||\cdot||_{L
\left(N_{\psi} \right)} =  ||\cdot||_{L (N)}$. It was proved that
there exist $ \infty > C_2 = C_2(\psi) \ge C_1 = C_1(\psi) >0 $ such that
for arbitrary r.v. $ \ \xi \ $

\begin{equation}\label{equivalence norm}
 C_1 ||\xi||_{G\psi} \le||\xi||_{L (N)}  \le C_2 ||\xi||_{G\psi}.
\end{equation}

 \ Of  course, the last relation  has meaning iff $ \ ||\xi||_{L (N)} < \infty  \  $ or equally $ \ ||\xi||_{G\psi} < \infty. \ $

\vspace{4mm}

\begin{example}
{\rm If for instance $  \ \psi(p) = \psi_m(p)\stackrel{def}{=}
p^{1/m}, \ p \in [1, \infty)$, \ $m = \const > 0$, then

$$
 0 \ne \xi \in G\psi_m \Leftrightarrow \ T_{\xi}(u) \le \exp \left(-C(m) u^m \right), \ u \ge 1.
$$
Define also the correspondent Young - Orlicz function
$$
N_m(u) := \exp \left( |u|^m \right), \  |u| \ge 1; \ \  N_m(u) = e  \cdot u^2, \ |u|< 1.
$$
 The relation \eqref{equivalence norm} means in addition in this case
\begin{equation}\label{equivalence norm example}
||\xi||G\psi_m \le C_1(m) ||\xi||L(N_m)  \le C_2(m)
||\xi||G\psi_m, \ 0 < C_1(m) < C_2(m) < \infty.
\end{equation}
}
\end{example}

 \vspace{4mm}
Let us define the following $ \Phi$-function

$$
\phi_m(\lambda) =  |\lambda|^m, \ |\lambda| \ge 1; \  \phi_m(\lambda) = \lambda^2, \ |\lambda| < 1.
$$

The Orlicz norm is quite equivalent on the set of mean zero random
variables to the $ \ B(\phi_m) \ $ one,   but only in the case when
$ \ m \ge 1. \ $
  \ Notice that  in the case when $  \ m \in (0,1) \ $ the correspondent random variable $ \ \xi \ $ does not satisfy in general case the Cramer's condition.   Therefore,
it can not belongs to arbitrary $ \ B(\phi) \ $    space. \par

\vspace{3mm}

 \section{Main result: preliminary upper estimate. }

\vspace{3mm}

 Let $\phi (\cdot) \in \Phi$  and let $\{ \xi_n \}, \ n = 1,2,\ldots $ be an arbitrary sequence of normed r.v. from the
space $B(\phi)$ such that

\begin{equation}\label{Bnorm1}
||\xi_n||_{B(\phi)} = 1.
\end{equation}

 \ Denote $ \ r_n := \nu^{-1}( n), \ u_0 := \nu^{-1}(1), \ \theta_n:= \xi_n/r_n. \ $   Further, let us impose the following condition of super - multiplicativity
 on the function $ \ \nu = \nu(u): \ $

\begin{equation}\label{condition nu super-multiplicative}
\exists u_1 = \const > 0 \  \ : \ \ \forall a,b \ge u_1 \
\Rightarrow \nu(ab) \ge  \nu(a) \ \nu(b).
\end{equation}

 \ We take as the value $ \ u_1 \ $ its minimal non - negative one. \par
 \  Define also by $ \ k_1 = k_1(u, \nu(\cdot)) \ $ the (fixed) minimal positive integer number greatest or equal than $ \nu(u_1): \ $

\begin{equation}\label{kEquation}
    \ k_1 = [\nu(u_1)] = k_1(\nu, \ u_1)  \stackrel{def}{=} \min \{ \ m, \ m = 1,2,3,\ldots: \ m \ge \nu(u_1) \  \},
\end{equation}

 and put also  $ \ u_2 := \max(u_0, \ u_1). \ $ \par

 \ Define also the following tail functions

\begin{subequations}
\begin{equation}\label{Pn}
 P_n (u):= {\bf P}(\theta_n > u),  \  P(u):= \sup_{n \ge k_1}
P_n(u) = \sup_{n \ge k_1} {\bf P}(\theta_n > u),
\end{equation}

\begin{equation}\label{Pn segnato}
\overline{P}_n(u) := {\bf P}(\max_{i=k_1, k_1 +1,k_1+2,\ldots, k_1+n} \theta_i
> u), \ \overline{P}(u) := \sup_{n \ge k_1}  \overline{P}_n(u),
\end{equation}
\begin{equation}\label{Pn+}
P^+_n(u) := {\bf P} \left( \  \max_{i=k_1}^{k_1 +n} \left[ \ \frac{\xi_i}{\nu^{-1}(n)} \ \right] > u  \  \right), \ P^+(u) := \sup_{n \ge k_1} P^+_n(u).
\end{equation}
\end{subequations}

\vspace{5mm}

 \begin{theorem}
  \ Suppose that the function  $ \ \nu(\cdot) = \phi^*(\cdot) \ $
  be the  Young - Orlicz function defined in \eqref{Young-Fenchel transform} and satisfies
 the condition of super - multiplicativity \eqref{condition nu super-multiplicative}.\par
  \ Suppose also that the sequence of random variables $ \ \{  \xi_n \} \ $
      satisfies  the norming condition \eqref{Bnorm1}.  We state that

\begin{subequations}
\begin{equation}
P_n(u)  \le \exp \left( \ - n \ \nu(u) \  \right),  \ \ P(u) \le \exp \left( \ - k_1 \nu(u) \  \right), \ n \ge k_1, \ u \ge u_1;
\end{equation}

\begin{equation}
\overline{P}(u) \le (1- 1/e)^{-1}  \exp \left( \ -  k_1 \ \nu(u) \  \right), \ u > u_2;
\end{equation}

\begin{equation}
P_n^+(u) \le  n \ \exp ( - n \ \ \nu(u) ), \
P^+(u)  \le \exp \left( \ -  \nu(u) \  \right), \ u > u_2.
\end{equation}
\end{subequations}
\end{theorem}
\noindent

\begin{proof} Let us  investigate at first the probability $  \ P_n = P_n(u) \ $. We have using the conditions \eqref{condition nu
super-multiplicative},  (3.2) for all the sufficiently large values $ \ u \ge u_1 \  $ and $ \ n \ge k_1 \ $

$$
P_n(u) = {\bf P}(\xi_n/r_n > u)\le \exp \left( \ - \nu(u \ \nu^{-1}(n)) \   \right) \le
$$
$$
\exp \left( \ - \nu(u) \ \nu(\nu^{-1}(n))  \ \right) = \exp \left( \ - \nu(u) \ n  \ \right).
$$

 \ Further,

$$
\overline{P}(u) = {\bf P}  \left( \ \cup_{j=k_1}^{\infty} \{ \xi_j/r_j > u   \}  \ \right) \le
$$

$$
\sum_{j=k_1}^{\infty} {\bf P}  \left( \  \xi_j/r_j > u  \ \right) \le
\sum_{j=k_1}^{\infty} \exp \left( \ - \nu(u \ r_j) \  \right) \le
$$

$$
\sum_{j=k_1}^{\infty} \exp \left( \ - \nu(u) \ j  \ \right) = \left(1 - e^{-\nu(u)} \right)^{-1} \ \exp \left( - k_1 \ \nu(u) \  \right) \le
$$

$$
 (1-1/e)^{-1} \exp \left( \ -  k_1 \ \nu(u) \  \right), \ u > u_2.
$$

 \ Let us estimate now the probability $ \ P_n^+(u). \ $ We deduce acting analogously

 $$
 P_n^+(u) = {\bf P} \left( \ \frac{\max_{i=k_1}^{k_1 +n} \ \xi_i}{\nu^{-1}(n)} > u \ \right) = {\bf P} \left( \ \max_{i=k_1}^{k_1 +n} \xi_i > u \ \nu^{-1}(n) \ \right) =
 $$

$$
 {\bf P} \left( \  \cup_{i=k_1}^{k_1 +n}  \{ \xi_i > u \ \nu^{-1}(n) \} \ \right) \le  \sum_{i =k_1}^{k_1 +n}  {\bf P} \left( \ \xi_i >  u \ \nu^{-1}(n)  \ \right)  \le
$$

$$
\sum_{i=k_1}^{k_1 + n} \exp \left(  - \nu(u \ \nu^{-1}(n) \right) \le \sum_{i=k_1}^{k_1 + n} \exp ( - n \ \nu(u)  ) =
$$

$$
n \ \exp ( - n \ \ \nu(u) );
$$
therefore
$$
P^+(u)  \le \exp(- \nu(u)), \ \  u \ge u_2,
$$
as long as $ \ \nu(u) \ge 1. \ $ \par

\end{proof}

\vspace{4mm}

 \section{Main result: more fine upper estimates. }

\vspace{5mm}

 \  {\it We assume in this section only the exponential inequality of the form}

 \vspace{3mm}

\begin{equation}\label{exponential inequality}
{\bf P}(\xi_i > x) \le \exp ( \ - \nu(x) \ ), \ \ x \ge 1,
\end{equation}
{\it with suitable (convex) Young - Orlicz function $ \  \nu  = \nu(x)  \  $ having
 continuous differentiable  strictly increasing to infinity  derivative function} $ \ \nu'(x). \ $ \par
 \ We emphasize that the r.v. $ \ \{ \xi_i  \}, \ i = 1,2,3,\ldots \ $   are  $ "ad \ lib" \ $ dependent. \par
 \ Denote also for brevity $ \  \vec{\xi} \stackrel{def}{=} \{ \xi_1, \xi_2, \xi_2,  \ldots  \}.  \  $ \par

\vspace{3mm}

 \ Let us define  the following variables
\begin{equation}\label{vn-wn}
r_n = \nu^{-1}( n),  \ \ \ w_n = \frac{1}{\nu'(r_n)} =
\frac{1}{\nu'(\nu^{-1}(n))},  \ \  \ n \ge 3,
\end{equation}

\begin{equation}\label{xin-ron}
\overline{\xi}_n = \max_{i=1,2,\ldots,n}\xi_i, \ \  \ \rho_n =
\frac{\overline{\xi}_n - r_n}{w_n} .
\end{equation}

 \ The variables $ \ \{ \rho_n \} \ $ are the sequences of random variables (r.v.).
 Let us estimate the {\it uniform}  tail function for ones.  \par

\begin{theorem}\label{th estimate tail function rhon}
 \ Let $\{\xi_i\}$ be a sequence of random variables
satisfying the condition \eqref{exponential inequality} and let
$\rho_n$ be defined in \eqref{xin-ron}. Then the
r.v. $\overline{\xi}_n = \max_{i=1}^n \xi_i$
has the following representation
\begin{subequations}\label{representation rhon}
\begin{equation}\tag{\ref{representation rhon}}
\overline{\xi}_n  = \nu^{-1}(\ln n) +
\frac{\rho_n}{\nu'(\nu^{-1}(\ln n))},  \ \ \ n \ge 3,
\end{equation}
 wherein
\begin{equation}\label{estimate tail function rhon}
\sup_{n \ge 3} {\bf P}(\rho_n > u) \le e^{-u}, \ \ \ u \ge 0.
\end{equation}
\end{subequations}
\end{theorem}

\noindent \begin{proof} The representation \eqref{representation rhon}
follows  from the direct definition of $\rho_n$. Further,

\begin{equation}
\begin{split}
{\bf P}(\rho_n > u) & =  {\bf P} \left(   \frac{\overline{\xi}_n -
r_n}{w_n} > u \right)  = {\bf P} \left( \overline{\xi}_n > r_n + u \
w_n   \right)\\
& ={\bf P}(\cup_{i = 1,2,\ldots, n} \{ \xi_i > r_n + u \ w_n \} )
\le \sum_{i=1}^n {\bf P}( \xi_i > r_n + u \ w_n  )\\ \\
&\leq n \exp
( - \nu(r_n + u \ w_n)  )  =  \exp \left( \ln n -
\nu(r_n + u \ w_n)   \right)\\ \\
&\leq \exp( \ln n - \nu(r_n) - \nu'(r_n) \ w_n \ u)  = \exp(-u),
\end{split}
\end{equation}
from which we get \eqref{estimate tail function rhon}. \par

\end{proof}

\vspace{4mm}

 \ Introduce the (deterministic) sequence
\begin{equation}\label{def zn}
 z_n := \nu^{-1}(n) \cdot \nu'(\nu^{-1}(n))=\frac{r_n}{w_n},
\end{equation}
so that

\begin{equation}
\frac{\overline{\xi}_n}{\nu^{-1}(\ln n) }  = 1  +
\frac{\rho_n}{z_n},  \ \ \ n \ge 3,
\end{equation}
and define also the deterministic variable $ \ K = K[ \vec{\xi},
\nu] := \ $ as

\begin{equation}\label{variable K}
K = K[ \vec{\xi}, \nu] :=\inf \left\{ \ Y  > 0: \  \forall \epsilon
> 0 \ \Rightarrow \ \sum_{n=1}^{\infty} \exp \left[ \ - \ (Y +
\epsilon) \ z_n \ \right] < \infty \
  \right\};
\end{equation}
the case when $ K[ \vec{\xi}, \nu] = 0$ is not excluded and will be investigated further. \par

 \ The following Theorem is an immediate consequence of the well-known
lemma of Borel - Cantelli. \par

\begin{theorem}\label{th lim sup prob one}
 \  Let $\{\xi_i\}$ a sequence of random variables satisfying the condition \eqref{exponential inequality} and let as above $\overline{\xi}_n = \max_{i=1}^n \xi_i$.
 Let also the "constant"  $K$ be defined by \eqref{variable K}. Then

\begin{equation}\label{probability  limsup K}
{\bf P} \left( \ \overline{\lim}_{n \to \infty} \  \left[ \
\frac{\overline{\xi}_n}{\nu^{-1}(\ln n) } \ \right] \le 1 + K[
\vec{\xi}, \nu] \  \right) = 1.
\end{equation}
\end{theorem}
\vspace{5mm}

 \ Let us consider a more general case, namely, when the r.v. $\xi_i$
 are not necessarily identically distributed:
\begin{equation}\label{var not necessarily identically distributed}
{\bf P}(\xi_i > x) \le \exp ( \ - \nu_i(x) \ ), \ \ \  x \ge 1
\end{equation}
with continuous differentiable convex having  strictly increasing to
infinity  functions $ \ \nu_i'(x). \ $ \par

\vspace{4mm}

 \ Let us introduce now a modified  notation. Define the value $ q_n$
as the unique positive root of the equation
\begin{equation}\label{equation define vn}
\sum_{i=1}^n e^{ \ - \nu_i(q_n)   \ } = 1.
\end{equation}

 \ Note that in the case when $ \ \nu_i(\cdot), \ i = 1,2,\ldots  \ $ are equal $ \ \nu_i = \nu, \ $ the value $ \ q_n \ $
coincides with  introduced before value $ \ r_n. \ $ \par

\vspace{3mm}

 \ Put similarly as in \eqref{vn-wn}, \eqref{xin-ron} and \eqref{def zn}
\begin{equation}\label{wn-rhon-zn independent}
w_n := \frac{1}{\sum_{i=1}^n \nu'_i(q_n)}, \ \ \ \  \rho_n =
\frac{\overline{\xi}_n - q_n}{w_n} , \ \ \ \ z_n := \frac{q_n}{w_n}.
\end{equation}
and define the variable $K = K[ \vec{\xi}, \nu]$ as in
\eqref{variable K}. We get analogous results to Theorems \ref{th
estimate tail function rhon} and \ref{th lim sup prob one}.\par

\begin{thma}
{\it Let $\{\xi_i\}$ be a sequence of  random variables
satisfying the condition \eqref{var not necessarily identically
distributed}. Let $q_n$ be defined by \eqref{equation define vn} and
$w_n$ and  $\rho_n$  defined in \eqref{wn-rhon-zn independent}. Then
the r.v. $ \ \overline{\xi}_n = \max_{i=1}^n \xi_i \ $ has the
following representation}
\begin{subequations}\label{representation rhon independent}
\begin{equation}\tag{\ref{representation rhon independent}}
\overline{\xi}_n  = q_n +  w_n \ \rho_n,  \ \ \ n \ge 1,
\end{equation}
{\it and}
\begin{equation}\label{estimate tail function rhon independent}
\sup_{n \ge 3} {\bf P}(\rho_n > u) \le e^{-u}, \ u > 0.
\end{equation}
\end{subequations}
\end{thma}

\begin{thma2}
{\it Let $\{\xi_i\}$  be a sequence of independent random variables
satisfying the condition \eqref{var not necessarily identically
distributed} and let $q_n$ be defined by \eqref{equation define vn}. Then}
\begin{equation}\label{prob th 4-2a}
{\bf P} \left( \ \overline{\lim}_{n \to \infty} \  \left[ \
\frac{\overline{\xi}_n}{q_n } \ \right] \le 1 + K[ \vec{\xi}, \nu] \
\right) = 1. \end{equation}
\end{thma2}

\vspace{4mm}

\begin{remark}
 {\rm As far as we know, the statements of Theorems 4.1 - 4.2a are known for Gaussian variables, see \cite{Belyaev-Piterbarg1978}, \cite{Pickands1967}, \cite{Sgibnev 1996},
 \cite{Stout 1974}, \cite{Sung 1996}, \cite{Welsch 1973} etc.} \par
\end{remark}

\vspace{4mm}

\begin{remark}
{\rm We investigate separately the possible case when $ \ K = 0, \ $ i.e. when

\begin{equation}\label{prob th 4-3b}
{\bf P} \left( \ \overline{\lim}_{n \to \infty} \  \left[ \
\frac{\overline{\xi}_n}{q_n } \ \right] \le 1  \ \right) = 1.
\end{equation}

 \ The sufficient condition for this conclusion is the following:

\begin{equation} \label{Kzero}
\forall \epsilon > 0 \ \Rightarrow \sum_{n=3}^{\infty} e^{-\epsilon \ z_n} < \infty.
\end{equation}

  \ In turn, the last condition  is satisfied if for example

 \begin{equation}
\nu(x) = \nu_{(s)}(x) \stackrel{def}{=} \exp \left( C \ |x|^s \right) - 1, \ x \in R,  \ C,s = \const > 0,
 \end{equation}
or more generally when

 \begin{equation}
\nu(x) = \nu_{(s_1, s_2)}(x) \stackrel{def}{=} \exp \left( C_{1,2} \ |x|^{s_1} \ \ln^{s_2}(x) \right) - 1, \  x \ge e,
 \end{equation}
 where $  \ C_{1,2}, \  s_1 = \const > 0, \ s_2 = \const \in R.   \ $ \par
}
\end{remark}

\vspace{4mm}

\begin{remark}
{\rm Let us  investigate the case when

\begin{equation}\label{subseqk}
{\bf P} \left( \ \overline{\lim}_{k \to \infty} \  \left[ \
\frac{\overline{\xi}_{n(k)}}{q_{n(k) }} \ \right] \le 1  \ \right) = 1
 \end{equation}
for  some deterministic integer strictly increasing sequence $ \ \{n(k)\}, \ k = 1,2,\ldots. \ $ \par
  \ The sufficient condition for this conclusion is follow:

\begin{equation} \label{Kzero}
\forall \epsilon > 0 \ \Rightarrow \sum_{k=1}^{\infty} e^{-\epsilon \ z_{n(k)}} < \infty.
\end{equation}
}
\end{remark}

\vspace{4mm}

\begin{remark}
{\rm  An interest open problem: find conditions (necessary conditions and sufficient ones) for the relation

\vspace{3mm}

\begin{equation}\label{exactval1}
{\bf P} \left( \ \overline{\lim}_{n \to \infty} \  \left[ \
\frac{\overline{\xi}_{n}}{q_{n}} \ \right] = 1  \ \right) = 1.
 \end{equation}

 \vspace{3mm}

 \ To make sure that the problem is not simple, let us bring  next example.  Let $ \ \xi \ $ be some fixed r.v.  such that

\begin{equation}  \label{exacttailfunk}
{\bf P} (\xi \ge x) = e^{-\nu(x)}, \ x \ge 0,
\end{equation}
where as before $ \ \nu(\cdot) \in \Phi. \ $ One can choose for instance $ \ \nu(x) = x^2/2, \ x \in R. \ $ \par
 \ Define the "sequence" $ \ \{ \xi_i\}, \ i = 1,2,3,\ldots \ $ of normed random variables  for which $ \ \xi_i = \xi \ $ for all the values $ \ i. \ $ \par
 \ We have here

\begin{equation}\label{exactvalue0}
{\bf P} \left( \ \overline{\lim}_{n \to \infty} \  \left[ \
\frac{\overline{\xi}_{n}}{q_{n}} \ \right] = 0  \ \right) = 1.
 \end{equation}

 \vspace{3mm}

}
\end{remark}

\vspace{3mm}

 \ Let us consider a more complicated problem: under which conditions on the sequence of r.v. $ \ \{\xi_i\}, \ i = 1,2,\ldots \ $ the mentioned before upper limit
is greatest or equal 1? \par

\vspace{3mm}

 \ We introduce the following condition of {\it  supermultiplicativity  } on the function $ \ \nu = \nu(x): \ $%
\begin{equation}  \label{supermultiplicativity}
 \exists u_4 \in (1,\infty), \ \exists \epsilon_1 \in (0,1) \ :  \ \forall \epsilon \in (0,1) \, ,  \ \forall A \ge u_4  \
\Rightarrow \ \nu(A \cdot ( 1 - \epsilon)) \le \nu(A) \cdot \nu(1 -
\epsilon_1)
 \end{equation}

\vspace{3mm}

 \ Recall that the sequence $ \ \{q_n\} \ $ is in the sequel defined in \eqref{equation define vn}. \par

\vspace{5mm}

\begin{theorem}\label{lim ge one}
 \ Let $\{\xi_i\}$ be a sequence of {\it independent} random variables
satisfying the condition \eqref{exacttailfunk} such that the
correspondent function $ \ \nu = \nu(x) \ $ satisfies the condition
\eqref{supermultiplicativity}. Then

\begin{equation}\label{exactvalue0}
{\bf P} \left( \ \overline{\lim}_{n \to \infty} \  \left[ \
\frac{\overline{\xi}_{n}}{q_{n}} \ \right] \ge 1  \ \right) = 1.
 \end{equation}
\end{theorem}

\vspace{4mm}

 \ {\bf Proof.} It is sufficient  by virtue of independence to ground that

\begin{equation} \label{Sigmaeps0}
\forall \epsilon \in (0,1) \ \Rightarrow \Sigma(\epsilon) = \infty,
\end{equation}
 where

\begin{equation}\label{Sigmaepsilon1}
\Sigma(\epsilon) \stackrel{def}{=} \sum_{n=3}^{\infty} {\bf P}(\xi_n/q_n > 1 - \epsilon).
\end{equation}

 \ Let us estimate from below the value $ \ \Sigma(\epsilon), \ \epsilon \in (0,1) \ $ from the relation \ref{Sigmaepsilon1}. We have taking into account the
 condition \eqref{supermultiplicativity}:

\begin{equation}
\Sigma(\epsilon) = \sum_{n=3}^{\infty} {\bf P} \left( \ \xi_n > [ q_n \cdot ( 1 - \epsilon) ] \ \right) = \sum_{n=3}^{\infty} \exp \left( - \nu(\nu^{-1}(\ln n) \cdot (1 - \epsilon) ) \right) \ge
\end{equation}

\begin{equation}
\sum_{n=3}^{\infty} \exp(- \ln n \cdot (1 - \epsilon_1)) = \sum_{n=3}^{\infty}  n^{-(1 - \epsilon_1)} = \infty,
\end{equation}
Q.E.D.

\vspace{4mm}

\begin{example}
 {\rm Suppose in  addition that $ \ \nu(x) = m^{-1} \ x^m,  \ \ x \ge 1$, \ \  $ m > 1, \ $ i.e.
\begin{equation}
{\bf P}(\xi_i > x) \ \le \exp \left( \ - x^m/m \   \right), \ x \ge
1.
\end{equation}
 We deduce after simple calculations for the values $ \ n \ge 3 \ $
\begin{equation}
q_n = (m \ \ln n)^{1/m}, \ \  w_n = (m \ \ln n)^{1/m - 1}, \ \
z_n=\frac{q_n}{w_n} = m \ \ln n,
\end{equation}
so that, by \eqref{representation rhon independent},
\begin{subequations}\label{xin example4-1}
\begin{equation}\tag{\ref{xin example4-1}}
\overline{\xi}_n = (m \ \ln n)^{1/m} + \frac{\rho_n}{ (m \ \ln n)^{1
- 1/m}},
\end{equation}
where
\begin{equation}\label{prob example 4-1}
\sup_{n \ge 3} {\bf P}(\rho_n > u) \le e^{-u}, \ u > 0,
\end{equation}
\end{subequations}
and by \eqref{prob th 4-2a}, with probability one, we have
\begin{equation}
\overline{\lim}_{n \to \infty} \frac{\overline{\xi}_n}{ (m \ln
n)^{1/m}  } \le 1 + \frac{1}{m}.
\end{equation}
For the independent centered (mean zero) r.v. $ \ \xi_i, \ i =
1,2,\ldots \ $ with distribution
\begin{equation}
{\bf P}(\xi_i > x)  = \exp \left( \ - x^m/m \   \right), \ x \ge 1
\end{equation}
one can deduce
\begin{equation}
\overline{\lim}_{n \to \infty} \frac{\overline{\xi}_n}{ (m \ln
n)^{1/m}  } = 1.
\end{equation}}
\end{example}

\vspace{4mm}

 \ Further, introduce the following deterministic increasing sequence

$$
n_0(k) \stackrel{def}{=} \Ent \left[ k^{\Delta(k)} \right], \ k = 1,2,\ldots,
$$
where $ \ \Ent(Z) \ $ denotes the integer part of a positive number
$ \ Z, \ $ and $ \ \{\Delta(k)\}, \ k = 1,2,\ldots \ $ is arbitrary
positive  non - random strictly increasing to infinity numerical
sequence: $ \ \Delta(k+1) > \Delta(k), \ \lim_{k \to \infty}
\Delta(k) = \infty. \ $ We deduce for the considered here random
variables with probability one

$$
\overline{\lim}_{k \to \infty} \frac{\overline{\xi}_{n_0(k)}}{ (m \ln n_0(k))^{1/m}  } \le 1.
$$

\vspace{5mm}

\begin{example}\label{example subgaussian}
{\rm If all the i.d. random variables $ \ \xi_i \ $ are in
additional subgaussian ($m = 2$ ), then evidently
\begin{equation}
\overline{\xi}_n = (2 \ \ln n)^{1/2} + \frac{\rho_n}{ (2 \ \ln
n)^{1/2}},
\end{equation}
where as before
\begin{equation}
{\bf P}(\rho_n > u) \le e^{-u}, \ u \ge 1,
\end{equation}
and with probability one
$$
\overline{\lim}_{n \to \infty} \frac{\overline{\xi}_n}{ (2 \ln n)^{1/2}  } \le 3/2.
$$
If in addition the r.v. $ \ \xi_i \ $ are independent and
$$
{\bf P}(\xi_i > x) = \exp \left(  \  - x^2/2 \ \right), \ x \ge 1,
$$
 then
$$
\overline{\lim}_{n \to \infty} \frac{\overline{\xi}_n}{ \ (2 \ln n)^{1/2} \ } = 1
$$
almost everywhere.}
\end{example}

\vspace{5mm}

\begin{remark}
{\rm The condition \eqref{var not necessarily identically
distributed} is satisfied  in the  following important case:
\begin{equation}\label{prob remark 4-2}
{\bf P}( \xi_i > x) \le C_1 \ i^{\kappa} \ e^{-\nu(x)},  \ \ \ x \ge 1,
\end{equation}
where $ \ C_1 \in (0, \ \infty), \ \kappa = \const \ge 0 \ $ and the function $ \ \nu(\cdot) \ $ is described before.}
\end{remark}
This case take place  in turn in the theory of random fields.
Namely, let $ \ Z \ $ be an arbitrary set, for instance, $ \ Z =
R^d_+, \ $  and let $ \ \{ T_i \}, \ i = 1,2,\ldots \ $
 be an increasing complete {\it sequence} of subsets of $ \ Z: \ $
 \begin{subequations}\label{sequence sets}
\begin{equation}\tag{\ref{sequence sets}}
T_1 \ne \emptyset, \ T_i \subset  T_{i+1}, \ \cup_{i=1}^{\infty} T_i
= Z.
\end{equation}
 Let also $ \ \zeta(z), \ z \in  Z \ $   be a separable numerical valued random field (process).  Put
\begin{equation}\label{separable r v}
\xi_i := \sup_{z \in T_i} \zeta(z).
\end{equation}
\end{subequations}
 \ The estimation of the form \eqref{prob remark 4-2} is obtained in particular by
means of the modern method of majorizing measures under appropriate
natural conditions, see in particular \cite{Ostrovsky-Sirota Oct2013} .\par

 \vspace{3mm}

 \ So, assume the estimate \eqref{prob remark 4-2} be given. We find
consequentially as $ \  n \to \infty \ $
\begin{eqnarray*}
 q_n \sim v^{-1}( C_2(\kappa) + (\kappa + 1) \ \ln n  ), \ \ \  w_n = \frac{1}{n \ \nu'(r_n)}.
\end{eqnarray*}
 \ It remains to apply the proposition of Theorem 4.1a. \par

 \vspace{3mm}

\noindent If in addition $ \ \nu(x) = \nu_m(x) := m^{-1} x^m, \ x
\ge 1, \ m > 1,   \ $ then
\begin{subequations}
\begin{equation}
q_n = r^{(m)}(n)  \sim  \left[ \ m( C_3 + (\kappa + 1) \ \ln n) \ \right]^{1/m},
\end{equation}
\begin{equation}
w_n = w^{(m)}(n) \sim m^{-1/m} \ n^{-1} \ \left[ \ C_3  + (\kappa +
1) \ \ln n \ \right]^{-1/m}.
\end{equation}
\end{subequations}

\vspace{4mm}

\noindent As a slight consequence under these conditions we have
\begin{equation}
\overline{\lim}_{n \to \infty} \left\{ \
\frac{\overline{\xi_n}}{v^{(m)}_n} \ \right\} \le C_4(m) = \const
\in (0, \infty),
\end{equation}
and the last estimate is essentially in general case non -
improvable. \par
The case when
\begin{equation}
{\bf P}(\xi_i > x) \le C_5 \  i^{\gamma} \ \exp ( \ - \nu_i(x) \ ),
\ \ \  x \ge 1, \ \ \gamma > 0
\end{equation}
may be investigated quite analogously. \par

 \vspace{3mm}

 \ The {\it lower bound} for the distribution of  the sequence of r.v. $ \ \rho_n \ $ under appropriate conditions is given in particular
 in the next section, see, e.g. \eqref{Lower bound}.  \par

\vspace{5mm}

 \section{Main result: lower estimates. }

\vspace{4mm}

 \ Let us show in this section an {\it unimprovability} in general case
of the obtained estimates. We consider for this purpose  the
sequence of {\it  independent} random variables $ \ \xi_i, \ i =
1,2, \ldots, \ $ with the following tail behavior

\begin{equation}\label{prob sec5-1}
{\bf P}(\xi_i > x) = \exp ( - \nu(x)), \ \ \ x \ge 1,
\end{equation}
where as before $ \  \nu(\cdot) \ $ is certain Young - Orlicz non - negative  twice
continuous differentiable  convex function such that its derivative
$ \ \nu'(x), \ $ as well as  itself $\nu(x)$, are strictly
increasing to infinity:

\begin{equation*}
\lim_{x \to \infty} \nu(x) = \lim_{x \to \infty} \nu'(x) = \infty,
\end{equation*}

\begin{equation*}
\nu(0) = 0, \ x > 0 \ \Rightarrow \nu'(x) > 0.
\end{equation*}

 \ The set of all such  functions will be denoted by $ \ \Phi: \ \Phi = \{ \nu(\cdot)  \}. \ $ \par

 \ Let $ \ \gamma_n \ $ be an arbitrary positive numerical bounded sequence tending to zero:

$$
  \lim_{n \to \infty}\gamma_n = 0, \ 0 < \gamma_n \le 1.
$$

 \noindent Let us introduce the  variables $ \ \epsilon_n, \ \Theta_n, \ R_n \ $
 from the following system of equations

\begin{subequations} \label{systemeq}
\begin{equation}
2 R_n := \sup_{u \in [0, \Theta_n]} | \ \nu^{''}(q_n + u \ w_n) \ |,
\end{equation}

\begin{equation}
\epsilon_n := R_n \ w^2_n, \ \Theta_n := \frac{\gamma_n}{\sqrt{\epsilon_n}},
\end{equation}
\end{subequations}
and recall that as before

$$
q_n = \nu^{-1}(\ln n), \ \rho_n = \left( \overline{\xi}_n - q_n \right)/w_n, \ w_n = 1/\nu'(q_n).
$$

\vspace{4mm}

\begin{theorem}\label{lower estimate}

 \ {\it We suppose that there exists a solution  of the last system  such that}

 \begin{subequations}\label{lim zero}
\begin{equation}\tag{\ref{lim zero}}
 \lim_{n \to \infty} \epsilon_n = 0, \ \lim_{n \to \infty} \gamma_n = 0,
\end{equation}
{\it and such that}
\begin{equation}\label{lim teta infty}
\lim_{n \to \infty} \Theta_n = \infty.
\end{equation}
\end{subequations}
 {\it Let us restrict ourselves to the following  interval for the values} $u:$
\begin{equation}\label{range teta}
u \in [1, \ \Theta_n].
\end{equation}

\vspace{4mm}

 \ {\it Then}

\begin{equation} \label{our prop}
  {\bf P} ( \rho_n > u) \ge e^{-\gamma^2_n} \cdot e^{-u} - e^{-2u}, \ 1 \le u \le \Theta_n.
\end{equation}
\end{theorem}
\vspace{4mm}

{\bf Proof.}

 \ We have applying the well known Bonferroni inequality  and taking into
account the independence
$$
 {\bf P} ( \rho_n > u) = {\bf P}(\overline{\xi}_n > q_n + u \ w_n) = {\bf P} \left( \cup_{i=1}^n \{ \ \xi_i > q_n + u \ w_n \ \}     \right)  \ge
$$
$$
\sum_{i=1}^n {\bf P}(\xi_i  >  q_n + u \ w_n ) - \sum \sum_{i,j =
1,2,\ldots,n; i \ne j}  {\bf P}(\xi_i  >  q_n + u \ w_n ) \cdot {\bf
P}(\xi_j  >  q_n + u \ w_n ) =:
$$
\begin{equation}\label{Sigma1-Sigma2}
\Sigma_1 - \Sigma_2.
\end{equation}

\noindent  Let us first estimate the value $\Sigma_1$. We have by \eqref
{prob sec5-1} and taking into account the restriction \eqref{range teta}

$$
\Sigma_1 = n \exp (- \nu( q_n + u \ w_n) ) = \exp \left\{  \ln n  - \nu(q_n + u \ w_n)  \right\} \ge
$$

\begin{equation}
\exp \left( \ - u - \epsilon_n u^2 \ \right)  \ge \exp \left( -u - \gamma^2_n \right) =
\exp \left( - \gamma^2_n  \right) \cdot \exp(-u).
\end{equation}

 \ As for the second term  $ \ \Sigma_2 \ $ in \eqref{Sigma1-Sigma2}:

\begin{equation}
\Sigma_2 \le n^2 {\bf P}^2(\xi_i > q_n + u \ w_n) \le n^2 \left[ \ \exp(-\nu(q_n + u \ w_n ))  \right]^2  \le
\end{equation}

\begin{equation}
n^2 \ \left\{ \ \exp \left[ \ - \ln n - u \  \right] \ \right\}^2 =
e^{-2 u}, \ u \ge 1.
\end{equation}

 \ Thus, we deduce for all the values $ \ u \ $ mentioned in
\eqref{range teta} under our assumptions and condition \eqref{lim
zero}

\begin{equation} \label{Lower bound}
{\bf P}(\rho_n > u) \ge \exp \left(- \gamma^2_n \right) \cdot \exp(-u)
- \exp ( -2 u ), \  \ \  u \in [1, \ \Theta_n],
\end{equation}
Q.E.D. \par

 \vspace{4mm}

 \ Both the propositions of the last two sections
 can be rewritten under the conditions formulated in this
section as follows.

 \vspace{4mm}

  \ {\bf Corollary 5.1.} Suppose that the sequence of r.v. $ \  \{\xi_i\}, \ i = 1,2,3, \ldots \ $ satisfies the condition (5.1), where
 $ \nu(\cdot) \in   \Phi. $ Then

 \begin{equation} \label{Bilateral bound}
\overline{\lim}_{u \to \infty} \overline{\lim}_{n \to \infty} \sup_{\nu \in \Phi} \left[ \ e^u \ {\bf P}(\rho_n > u) \ \right] = 1.
\end{equation}

\vspace{3mm}

 \ Note that the conditions of theorem (5.1) are satisfied for the function $ \ \nu(x) = \nu_{m,r}(x) \ $ of the form

\begin{equation}
\nu(x) = \nu_{m,r}(x) \stackrel{def}{=}  x^m \ [\ln x]^r, \ \ \ m >0, \ r \in R, \ \ \ x \ge e.
\end{equation}

\vspace{4mm}

 \section{The case of arrays. }

\vspace{4mm}

 \ Let us return to the announced case in section \ref{intro} of arrays of centered independent random variables.
 Let $ \{\xi_{n,i} \}, \ i = 1,2,\ldots,n; \ n = 1,2,\ldots \ $ be an
array of independent random variables with $ {\bf E} \xi_{n,i} = 0$
and $ \  0 <  \sigma_{n,i}^2 := {\bf E} \xi^2_{n,i} < \infty$,
satisfying \eqref{limsup1}.

 \ Recall the notation

$$
S_n := \sum_{i=1}^n \xi_{n,i}. \eqno(6.0)
$$

 \ Suppose that every r.v. $\xi_{n,i}$ satisfies the Cramer's
condition, on the other words, $\xi_{n,i}$ belongs to some
$B(\phi[n,i])$ space, where  \ $\phi[n,i](\cdot) \in \Phi$:

\begin{equation}
{\bf E} \exp \left( \ \pm \ \lambda \ \xi_{n,i} \   \right) \le \exp \left( \  \phi[n,i](\lambda)  \ \right),
\end{equation}
see \eqref{momentfunction},
referring  also to the limitation \eqref{mean infinity}. Of course, one
can take as the function $ \ \phi[n,i](\lambda)  $  the natural function
for the correspondent r.v. $ \  \xi_{n,i}. \ $ \par

 \ Introduce a new function, more precisely, the sequence of ones, also belonging to the set $ \ \Phi: \ $
\begin{equation}\label{funct chi}
\chi_n(\lambda) := \sum_{i=1}^n  \phi[n,i](\lambda),
\end{equation}
then the r.v. $S_n $ belongs to the space $B(\chi_n)$ and has therein the
norm which is less than 1:
\begin{eqnarray*}
{\bf E} \exp( \ \lambda \ S_n) &= &\prod_{i=1}^n {\bf E}
\exp(\lambda \ \xi_{i,n}) \leq
\prod_{i=1}^n \exp (\phi[i,n](\lambda)) =\exp(\chi_n(\lambda)),
\end{eqnarray*}
therefore
\begin{subequations}\label{prob array}
\begin{equation}\tag{\ref{prob array}}
{\bf P} (S_n > u) \le \exp (- \kappa_n(u)), \ \ \  u \ge 0,
\end{equation}
where
\begin{equation}\label{kn prob array}
\kappa_n(u) := \chi_n^*(u) = \sup_{\lambda \in \Dom[\chi_n(\cdot)]} (\lambda \ u -
\chi_n(\lambda)).
\end{equation}
\end{subequations}

 \ It remains to apply Theorems \ref{th estimate tail function rhon},
\ref{th lim sup prob one}. Indeed, the sequence of the r.v.
$\overline{S}_n  := \max_{i=1}^n S_i$ allows the following
representation alike in the fourth section
\begin{subequations}\label{Sn bar}
\begin{equation}\tag{\ref{Sn bar}}
\overline{S}_n  = \kappa_n^{-1}(\ln n) +
\frac{\rho_n}{\kappa_n'(\kappa_n^{-1}(\ln n))}, \ \ \ n \ge 3,
\end{equation}
where as above
\begin{equation}\label{prob section 6}
\sup_{n \ge 3} {\bf P}(\rho_n > u) \le e^{-u}, \ \ \  u > 0.
\end{equation}
\end{subequations}
 Put now
$$
y(n) = \kappa_n^{-1}(\ln n) \cdot \kappa_n'(\kappa_n^{-1}(\ln n)),
$$
so that
\begin{equation}
\frac{\overline{S}_n}{\kappa_n^{-1}(\ln n) }  = 1  +
\frac{\rho_n}{y(n)},  \ n \ge 3,
\end{equation}
and define the (non - random) variable $ \  L = L[ \vec{\xi},  \{
\kappa(\cdot) \} ] := \ $ as
\begin{equation}\label{r v L}
L = L[ \vec{\xi},  \{ \kappa(\cdot) \} ] :=\inf \left\{ Y > 0 \ : \
\forall \epsilon > 0 \ \Rightarrow \ \sum_{n=1}^{\infty} \exp \left[
\ - \ (Y + \epsilon) \ y(n) \ \right] < \infty \  \right\}.
\end{equation}

 \ The next statement follows immediately, as before, again from the
well known lemma of Borel - Cantelli. \par

\vspace{5mm}

\begin{theorem}
\begin{equation}
{\bf P} \left( \ \overline{\lim}_{n \to \infty} \  \left[ \
\frac{\overline{S}_n}{\kappa_n^{-1}(\ln n) } \ \right] \le 1 + L \
\right) = 1.
\end{equation}
\end{theorem}

\vspace{5mm}

\begin{example}
{\rm Assume that all the centered i.d. random variables $ \xi_{n,i}
$ are subgaussian and independent and set
\begin{equation}
\beta_{n,i} := ||\xi_{n,i}||\Sub \in (0,\infty).
\end{equation}
Define
\begin{equation}
\overline{\beta}_n = \left[ \ \sum_{i=1}^n \beta^2_{n,i} \
\right]^{1/2}.
\end{equation}
If $ \ \beta_{n,i} = 1, \ i = 1,2,\ldots,n, \ $ then it is easily to
verify that all the conclusions of example \ref{example subgaussian}
remains true. Let us consider now the general case. Denote
\begin{eqnarray*}
\overline{S}_n := \max_{i = 1,2,\ldots,n} \sum_{j=1}^i \xi_{n,j}.
\end{eqnarray*}
We conclude
\begin{equation}\label{Sn general case}
\overline{S}_n = \overline{\beta}_n \left( \sqrt{2 \ \ln n} +
\frac{\rho_n}{\sqrt{2 \ \ln n} }   \right),
\end{equation}
or equally
\begin{equation}
\frac{\overline{S}_n}{\overline{\beta}_n \ \sqrt{2 \ \ln n}} =  1 +
\frac{\rho_n}{2 \ \ln n },
\end{equation}
where as before
\begin{equation}
\sup_{n \ge 3} {\bf P}(\rho_n > u) \le \ e^{-u}, \ u \ge 1.
\end{equation}
 As a consequence:
\begin{equation}\label{limsup Sn general case}
\overline{\lim}_{n \to \infty}
\frac{\overline{S}_n}{\overline{\beta}_n \ \sqrt{2 \ \ln n}} \le
\frac{3}{2}.
\end{equation}
If in addition the (independent) r.v. $ \ \xi_{n,i} \ $ are strictly
subgaussian, then one can take in the relations \eqref{Sn general
case}, \eqref{limsup Sn general case}
\begin{equation}
\overline{\beta}_n = \sqrt{\Var(S_n)}.
\end{equation}

\vspace{3mm}

\noindent Moreover, one can estimate the following tail probability
\begin{equation}
 Y(z) \stackrel{def}{=}
 {\bf P} \left( \ \sup_{n \ge 2} \left[ \ \frac{\overline{S_n}}{\overline{\beta_n} \ \sqrt{2 \ \ln n }}  \ \right] \ge 1 + z \  \right),  \ \ \ z \ge 1.
\end{equation}
 which  may be used in statistics and in the Monte - Carlo method. Indeed, we deduce subject to our  limitations

$$
Y(z) =  {\bf P} \left(  \ \cup_{n \ge 2} \left\{ \frac{\rho_n}{2 \ \ln n} \ge 1 + z  \  \right\}    \ \right) \le
$$

$$
\sum_{n=2}^{\infty} {\bf P} \left( \ \frac{\rho_n}{2 \ \ln n}  \ge 1 + z \  \right) \le
\sum_{n=2}^{\infty} n^{-2z} \le 2^{1 - 2z}. \eqno(6.16)
$$

}

\end{example}

\vspace{4mm}

\begin{remark}

 \ {\rm It is no hard to deduce the equalities of the form

\begin{equation}
{\bf P} \left( \ \overline{\lim}_{n \to \infty} \  \left[ \
\frac{\overline{S}_n}{\kappa_n^{-1}(\ln n) } \ \right] \le 1 \right) = 1.
\end{equation}
or moreover

\begin{equation}
{\bf P} \left( \ \overline{\lim}_{n \to \infty} \  \left[ \
\frac{\overline{S}_n}{\kappa_n^{-1}(\ln n) } \ \right] = 1 \right) = 1
\end{equation}
for the {\it arrays} of random variables, alike ones in fourth section. The lower bounds for tail of distribution of arrays sums
is a particular case for ones for the ordinary sums obtained in fifth section. \par

}
\end{remark}

\vspace{5mm}

 \section{ A Grand Lebesgue Spaces approach. }

\vspace{5mm}

 \ Let $ \ \xi \ $ be some numerical valued r.v.  from certain    $ \  G\psi \ $  space,
$ \ \psi \in U \Psi, \ \Dom [\psi] = [1,b), \ b = \const \in (1, \ \infty] \ $ and assume $ \ ||\xi||G\psi = 1. \ $   We deduce using
 estimate \eqref{estimateT}

\begin{equation}
  T_{\xi} ( y) \le \exp \left( \ - h_{\psi}^* (\ln  y ) \ \right), \ \ \ y \ge e,
\end{equation}
where as before

\begin{eqnarray*}
h(p) = h[\psi](p)  \stackrel{def}{=} p \ \ln \psi(p), \ \ \  1 \le p< b.
\end{eqnarray*}

 \ We obtained for Grand Lebesgue Spaces the analogous tail relation \eqref{exponential
inequality}; it remains to apply the results of Section 4.\par

\vspace{3mm}

\noindent The case when  the function  $ y \to \exp \left( \ -
h_{\psi}^* (\ln  y ) \ \right), \ y \ge e $, does not satisfy the
condition \eqref{var not necessarily identically distributed}, i.e.
when the r.v. $\xi_{n,i}$ have only \lq\lq power decreasing tail\rq\rq of
distribution, is investigated partially in \cite[pp.
44-48]{Ostrovsky1999}. The case of a very hard tail behavior for a
r.v.    $ \ \xi_{n,i} \ $ is considered in \cite{Braverman1991},
\cite{Braverman1994}, \cite{Ng-Tang-Yang}, \cite{Sung 1996},
\cite{Samko-Umarkhadzhiev-addendum}, \cite{Sgibnev 1996}, etc. \par

 \ Let us consider the following simple example. Suppose the r.v. $ \xi_{n,i}$ are such that
\begin{equation}\label{condition power tail}
\forall x \ge 1  \ \Rightarrow \   \sup_{n,i} {\bf P} \left( \
\xi_{n,i} > x \ \right) \le  x^{-p}, \ \  \ p> 0.
\end{equation}

 \ On the other words, $ \xi_{n,i}$ belong (uniformly in $ \ (n,i) \ )
\ $ to the unit ball of the so  - called Lorentz space $ \ L_{p,
\infty}$. Obviously, the condition \eqref{condition power tail} is
satisfied if
$$
\sup_{n,i} {\bf E} \left| \xi_{n,i}\right|^p  \le 1.
$$
 \ We deduce
\begin{subequations}\label{generalization Pizier1}
\begin{equation}\tag{\ref{generalization Pizier1}}
  {\bf P} \left( \ \overline{\xi}_n  > u \ n^{1/p} \ \right)  \le \sum_{i = 1}^n {\bf P} \left(\xi_{n,i} > u \ n^{1/p} \right) \le \\
     \sum_{i=1}^n \frac{1}{u^p \ n} = u^{-p}, \ \ \ u \ge 1,
\end{equation}
or equally
\begin{equation}\label{generalization Pizier2}
 \sup_{n = 1,2,3,\ldots} {\bf P} \left( \ n^{-1/p} \ \overline{\xi}_n  > u \ \right)  \le  u^{-p}, \ \ \ u \ge 1,
\end{equation}
\end{subequations}

 \ The last estimate is a slight generalization of one due by G. Pisier (\cite{Pisier 1980}). \par

\vspace{4mm}

\noindent Let's make sure that the estimates \eqref{generalization
Pizier1} and \eqref{generalization Pizier2} are essentially non -
improvable. One can choose for this purpose a {\it sequence} $\{
\xi_j\}, \ j = 1,2,\ldots$ \ of positive independent greatest than one identically
distributed random variables defined on suitable probability spaces
and such that
\begin{equation*}\label{lowbcond1}
{\bf P}(\xi_j \ge u) = u^{-p}, \ \ \  u \ge 1, \  \ \ p > 0.
\end{equation*}

\vspace{3mm}

  \ Let at first $ \ n = 1; \ $ then
\begin{equation*}
\sup_{n = 1,2,3,\ldots} {\bf P} \left( \ n^{-1/p} \ \overline{\xi}_n
> u \ \right)  \ge   {\bf P} (\xi_1 > u) = u^{-p}, \  \ \ u \ge 1.
\end{equation*}

 \ Let us consider now a general case $ \ n \ge 2. \ $ We have
consequently  applying once again the famous Bonferroni's inequality
\begin{equation*}
 {\bf P}(\overline{\xi}_n \ge u \ n^{1/p}) \ge  \sum_{j=1}^n {\bf P}(\xi_j \ge u) - \sum \sum_{i,j = 1,2,\ldots,n; \ i < j} {\bf P}(\xi_i \ge u, \ \xi_j \ge u) =
\end{equation*}
\begin{equation*} \label{lowbound1}
n \cdot \frac{1}{u^p \ n} - 0.5 n(n-1) \frac{1}{u^{2p} \ n^2}  =
u^{-p} - 0.5 \ u^{-2p}(1 - 1/n) \ge u^{-p} - u^{-2p},  \ \ \ u \ge 2.
\end{equation*}

\vspace{3mm}

 \ Let us consider now a more general case  when for some generating
function $ \ \psi \ \in \Psi \ \ \Rightarrow \xi_i \in G\psi, $ where $ \ \Dom[\psi] = [1,b), \ b = \const \in (1, \ \infty]; \ $
and assume  moreover
$$
\max_{i= 1,2,\ldots,n} ||\xi_i||_{G\psi} = 1.
$$
Define the functions
$$
g(u) := \ln \psi(1/y),\ \ \ y \in (1/b,1);
$$
$$
g_*(x) := \inf_{y \in (1/b,1)} (x  y + g(y)),
$$
which is named ordinary as
the so - called \lq\lq adjacent\rq\rq Young - Fenchel transform for the function $ \ g(\cdot).\ $ \par
 \ Obviously, the last function is correctly defined for all the values $ \  x \in {\bf R}. \ $  We have in particular, taking the value
 $ \ y_0 := (b+1)/(2b) \ $

\begin{equation}
g_*(x) \le (x y_0 + g(y_0)) = \left( x \ \frac{b+1}{2b} + g \left(  \frac{ b + 1}{2b} \right)   \right) < \infty.
\end{equation}

 \ But we need to use further only  positive values  for the variable $ \ x. \ $ \par

 \ It is proved in \cite[chapter1, section 1.10]{Ostrovsky1999}, that there exists a finite  "constant" $ \ \kappa_0(n) = \kappa_0[\psi](n) \ $ such that
   if $ \ \max_{i=2,3,\ldots,n} ||\xi_i|| \le 1, \ $ then

\begin{equation*}\label{majchar33}
|| \ \max_{i = 2, 3, \ldots,n} |\xi_i| \ ||_{G\psi} \le \kappa_0[\psi](n),
\end{equation*}
and herewith

\begin{equation*} \label{estimkapp0}
 \kappa_0[\psi](n) \le C(\psi) \ \exp(g_*(\ln n)), \ \ \ n \ge 2,
\end{equation*}
with correspondent  tail estimation \eqref{estimateT}. \

 \ The minimal value of the constant  $ \ \kappa_0[\psi](n), \ $  i.e. the value
\begin{equation*}\label{defkappa}
\kappa[\psi](n) = \kappa(n) \stackrel{def}{=}\sup_{\xi_i: \ \max( \
||\xi_i|| G\psi, \ i = 2,3,\ldots,n \ \in(0, \ \infty))} \left\{
\frac{||\max_{i = 2,3, \ldots,n} \ |\xi_i| \ ||G\psi}{\max_{i =
2,3,\ldots,n} \ ||\xi_i|| G\psi}   \right\}
 \end{equation*}
is named in a recent article \cite{Kozachenko at all 2018} as "M - characteristic" or "majorant characteristic"  for the space $ \ G\psi \ $ and alike spaces. The estimates of norm for maximum
$ \ || \ \max_{i = 2, 3, \ldots,n} |\xi_i| \ ||_{G\psi} $ common with suitable ones for $ \ \kappa(n), \ $
are used in \cite{Kozachenko at all 2018} as well as in the brochure  \cite{Ermakov etc. 1986} for the investigation of continuity for
random fields, conditions for Central Limit Theorem in the space of continuous functions and in turn in the parametric method Monte - Carlo. \par

\vspace{4mm}

 \ Let us investigate the  tail behavior for maximum distribution of (dependent, in general case)  random variables $ \ \xi_i, \ i = 1,2,\ldots \ $ from  certain Grand Lebesgue Spaces,
on the other hands,  having a heavy tails  of distribution. Namely, suppose

\begin{equation}
T_{\xi}(x) \le x^{-\alpha} \ L(x), \ \alpha = {\const} > 0, \ x \ge
1,
\end{equation}
where as before  $ \ L = L(x), \ x \ge 1 \ $ is some positive
continuous {\it slowly varying} function as $ x \to \infty $.
 Introduce an auxiliary function

\begin{equation}
M(y) = M_L(y) := \sup_{z \ge 1} \left\{ \ \frac{L(x \cdot z)}{L(z)} \  \right\},
\end{equation}

so that

\begin{equation}
L( x \ v ) \le L(v) \cdot M(x);
\end{equation}

then this function belongs also to the set $ SL: \ M_L(\cdot) \in SL. \ $ On the other words, this tail function is named as a regular varying ones.\par
 \ The correspondent $ \ \psi \ $ function is described in \cite{Kozachenko-Ostrovsky-Sirota Jan2017}, \cite{liflyandostrovskysirotaturkish2010}. \par

 \ Define the positive sequence $ \ U = U(n) \ $ so that $ \ U(1) = 1 \ $ and for the values$ \ n = 2,3,\ldots \ $  as a solution of an equation

\begin{equation}
U^{\alpha}(n) [L(U(n))]^{-1}  = n.
\end{equation}
 \ We deduce as before

$$
T \left[U^{-1}(n) \ \max_{i = 1,2,\ldots,n} \xi_i \right] (x) \le  \sum_{i=1}^n {\bf P}(\xi_i/U(n) > x) \le
n \ x^{-\alpha} U^{-\alpha}(n) \ L(U(n) \ x) \le
$$

\begin{equation}
n \ x^{-\alpha} \ U^{-\alpha}(n) \ L(U(n)) \ M(x) =  x^{-\alpha} \ M(x), \ x \ge 1.
\end{equation}

\vspace{3mm}

 \ To summarize: under the formulated above conditions

 \begin{equation}
\sup_{n \ge 1} T\left[U^{-1}(n) \ \max_{i = 1,2,\ldots,n} \xi_i \right] (x) \le  x^{-\alpha} \ M(x), \ x \ge 1.
 \end{equation}

\vspace{4mm}

 \ Further, let $ \ \upsilon = \upsilon(n) \ $ be any positive finite unbounded deterministic numerical sequence  for which

\begin{equation}
\sum_{n=2}^{\infty} \upsilon^{-\alpha}(n) \ M(\upsilon(n)) < \infty.
\end{equation}

 It follows immediately again from lemma of Borel - Cantelli that with probability one

\begin{equation}
\overline{\lim}_{n \to \infty} \left\{ \ \frac{\max_{i=1,2,\ldots,n} \xi_i}{U(n) \ \upsilon(n)} \ \right\} \le 1.
\end{equation}

\vspace{4mm}

 \ Moreover,

 \begin{equation}
 {\bf P} \left( \ \sup_{ n \ge 1} \frac{\overline{\xi}_n}{U(n) \ \upsilon(n)} \ge x \ \right) \le \sum_{n=1}^{\infty} {\bf P} \left( \frac{\overline{\xi}_n}{U(n) \ \upsilon(n)} \ge x \ \right) \le
 \end{equation}

\begin{equation}
x^{-\alpha} \ \sum_{n=1}^{\infty} \upsilon^{- \alpha}(n) \ M(x \cdot \upsilon(n)), \ x \ge 1.
\end{equation}

\vspace{4mm}

 \ Let us show now that our estimations obtained in this section are essentially non - improvable.  Consider the r.v. - s.  $ \ \xi_i, \ i = 1,2,\ldots \ $ such that

\begin{equation}
T_{\xi}(x) = x^{-\alpha} \ L(x), \ \alpha = {\const} > 0, \ x \ge 1,
\end{equation}
where as before  $ \ L = L(x), \ x \ge 1 \ $ is some positive continuous {\it slowly varying} as $ x \to \infty \ $ function. The following lower very simple estimate
holds true

$$
\sup_{n \ge 1} T\left[U^{-1}(n) \ \max_{i = 1,2,\ldots,n} \xi_i \right] (x) \ge  T\left[U^{-1}(1) \ \xi_1 \right] (x) =
$$

$$
T_{\xi_1} (x)  = {\bf P}(\xi_1 \ge x) = x^{-\alpha} \ L(x), \ x \ge
1;
 $$
herewith the case $ \  M(x) = L(x)  \ $ or at last when $ \  M(x) \le C \cdot L(x)  \ $ can not be excluded, for instance when

$$
L(x) = \ln^r(e \ x), \ r = \const > 0, \ x \ge 1
$$

 \vspace{4mm}

 \ Notice  that the case of {\it arrays} of the r.v.-s $ \ \xi_{i,n} \ $ under the same conditions  may be investigated quite analogously. \par

\vspace{4mm}

 \ Let us return to the source problem of estimation of partial sums for {\it independent  arrays} $ \ \{\xi_{i,n} \} \ $ of random variables, but now in the case of heavy
 tails of distributions:

\begin{equation}
S_n := \sum_{i=1}^n \xi_{i,n}, \  {\bf E} \xi_{i,n} = 0.
\end{equation}

  \ We assume that

\begin{equation}
\sup_n \max_{i = 1,2,\ldots,n} T_{\xi_{i,n}}(x) \le  T^{(\beta,
\gamma, L)}(x),  \beta = \const \in (2,\infty), \ \gamma = \const
\in \mathbf{R},
\end{equation}
where  (we recall)

\begin{equation}
T^{(\beta, \gamma, L)}(x) \stackrel{def}{=} x^{-\beta} \ (\ln
x)^{\gamma} L(\ln x), \ x \ge e.
\end{equation}
where as above  $ \ L = L(x), \ x \ge 1 \ $ is certain positive continuous {\it slowly varying} as $ x \to \infty \ $ function.  As we knew,

\begin{equation}
\sup_n \max_{i = 1,2,\ldots,n} ||\xi_{i,n}||G\psi^{(\beta,\gamma +1,L)} = C_1 = C_1(\beta,\gamma,L) < \infty.
\end{equation}

 \ One can apply the famous Rosenthal's inequality, see e.g. \cite{Naimark 2004}, taking into account the boundedness of correspondent Rosenthal's coefficient
 $ \ R(p) \ $ in the closed  segment $ \ p \in [1,\beta] \ $

\begin{equation}
\sup_n ||n^{-1/2} S_n||_p \le  C_2(\beta,\gamma,L) \ \psi^{(\beta,\gamma+1,L)}(p), \ 1 \le p < b,
\end{equation}
 or equally

\begin{equation}
 \sup_n ||n^{-1/2} S_n||G \psi^{(\beta,\gamma+1,L)} = C_3(\beta, \gamma,L) < \infty.
\end{equation}
 \ We conclude ultimately returning to the tail of distribution

$$
\sup_n T[n^{-1/2} S_n ](x) \le   C_4(\beta, \gamma,L) T^{(\beta, \gamma + 1, L)}(x)  =
$$

$$
 C_5(\beta, \gamma,L)  x^{-\beta} \ (\ln x)^{\gamma + 1} L(\ln x), \ x \ge e.
$$

\vspace{4mm}

 \ Further, let $ \ d = d(n) \ $ be a certain positive finite unbounded deterministic numerical sequence  for which  $ \ d(1) = 1 \ $ and

\begin{equation}
\sum_{n=2}^{\infty} \ d^{-\beta}(n) \ [\ln d(n)]^{\gamma +1} \ L(\ln d(n)) < \infty.
\end{equation}
  \ The last condition is in turn satisfied if for example

\begin{equation}
 d(n) \ge n^{1/\beta} \  [\ln n]^{\delta }, \ n \ge 2,
\end{equation}
 where $ \ \delta = \const  > (2 + \gamma)/\beta. \ $

\vspace{3mm}

 \  It follows immediately again from mentioned above lemma of Borel - Cantelli that with probability one

\begin{equation}
\overline{\lim}_{n \to \infty} \left\{ \ \frac{\max_{i=1,2,\ldots,n} S_n}{n^{1/2} \ d(n)} \ \right\} \le 1.
\end{equation}

\vspace{4mm}

 \ Moreover, we conclude as before

 $$
 {\bf P} \left( \ \sup_{ n \ge 1} \left\{ \ \frac{S_n}{n^{1/2} \ d(n) } \ \right\}  \ge x \ \right) \le \sum_{n=1}^{\infty} {\bf P} \left( \frac{S_n}{n^{1/2} \ d(n) \ } \ge x \ \right) \le
 $$

\begin{equation}
x^{-\beta} \ \sum_{n=1}^{\infty} \ d^{- \beta}(n) \ [ \ln(x \ d(n)) ]^{\gamma + 1} \ L (\ln (x \cdot \ d(n))), \ x \ge 1.
\end{equation}

 \vspace{4mm}

 \ Of course,   the norming function $ \ n^{1/2} \ d(n)  \  $  in the case of heavy tails of distribution of source  r.v. $ \ \xi_{i,n} \ $
significantly differs from the classical  ones $ \  n^{1/2} \ \ln \ln n  \ $ as well as $ \  n^{1/2} \ \ln n.  \ $ \par

 \ Note in addition that the {\it lower bound } for this probability is quite alike for obtained before. Indeed, assume that the r.v. - s $ \ \xi_{i,n}, \ $ not necessary independent
are such that  $ \ \forall x \ge e \ \Rightarrow \ $

\begin{equation}
 T_{\xi_{i,n}}(x) =  T^{(\beta, \gamma, L)}(x),  \beta = {\const} \in (2,\infty), \ \gamma = {\const} \in \mathbf{R},
\end{equation}
then

 \begin{equation}
 {\bf P} \left( \ \sup_{ n \ge 1} \left\{ \ \frac{S_n}{n^{1/2} \ d(n) } \ \right\}  \ge x \ \right) \ge {\bf P}(\xi_{1,1}  \ge x) = T^{(\beta, \gamma, L)}(x), \ x \ge e.
 \end{equation}

\vspace{5mm}

\section{Concluding remarks.}

\vspace{5mm}

{\bf A.} One can reduce the condition \eqref{condition nu
super-multiplicative} with a  more general one:
$$
\nu(a \ b) \ge R(a) \ \nu(b) \eqno(R)
$$
for some positive continuous  {\it increasing to infinity}  function $R = R(u)$ and for all sufficiently large values $a,b$. \par

\vspace{5mm}

\noindent {\bf B.} The condition \eqref{condition nu
super-multiplicative} is trivially satisfied if the function $ \nu(\cdot)$ has the form
$$
\nu(x) = \nu_{m,0}(x) = x^m, \ x \ge e,
$$
where we define
$$
\nu_{m,r}(x) \stackrel{def}{=}  x^m \ [\ln x]^r, \ \ \ m >0, \ r \in {\bf R,} \ \ \ x \ge e.
$$
In turn, the introduced function $ \nu_{m,r}(\cdot), \ r > 0$
satisfies the condition (R). \par

 \ In detail, it is no hard to calculate that when $ \ r > 0 \ $ then the function $ \ R(\cdot) \ $ may be choosed as

$$
R(x)  = C(r) \ x^m \ \ln^r(x), \ x \ge e,
$$
and $ \ R(x) := x^m, \ x \ge e, \ $ if $ \ r < 0. \ $ \par

 \ Let us bring again a more general example. Let $ \ L = L(x), \ x \ge 1 \ $ be some positive continuous {\it slowly varying} as $ x \to \infty \ $ function; the set of all
such a functions will be denoted by $ \ SV; \ SV = \{ L(\cdot) \}.  \ $  Recall that

$$
M(y) = M_L(y) := \sup_{z \ge 1} \left\{ \ \frac{L(x \cdot z)}{L(z)} \  \right\};
$$
then this function belongs also to the set $ SL: \ M_L(\cdot) \in SL. \ $ \par

 \ Define the following  Young - Orlicz  function $ \ \nu(\cdot) \ $

$$
\nu_{m,L}(x) := x^m \ L(x), \ x \ge 1;
$$
 then the correspondent $ \ R(\cdot) \ $ function from the condition (R) may be choosed in the form

$$
 R_{m,L} (y) = y^m \ M_L(y), \ y \ge 1.
$$

\vspace{5mm}

\noindent {\bf C.}  It is interest  by our opinion to investigate
also the case of {\it  continuous} "time", i.e. to describe the non
- asymptotical behavior  as $ \ T \to \infty \ $ of the random process

$$
\xi(T) := \sup_{t \in [0,T]} \zeta(t)
$$
or more generally

$$
\xi(T) := \sup_{t \in [0,T]} \int_0^t \zeta(s) \ \mu_T(ds),
$$
in the spirit of the classical LIL for Brownian motion. \par
\noindent Some preliminary results in this direction may be found in \cite[pp. 150 - 157]{Ostrovsky1999}. \\

\vspace{5mm}

\noindent {\bf D.} Let us show a possible application in statistics.
Consider for simplicity the following model:
 the r.v. $ \ \tau_i, \ i = 1,2,\ldots \ $ are i., i.d. r.v. with $ \ \theta =  {\bf E} \tau_i \in \mathbf{R} \ $ and
 $ \ \beta := ||\tau_i - \theta||\Sub  \in (0,\infty)$. Define the consistent ordinary estimate of the value $ \ \theta: \ $
$$
\theta_n := n^{-1} \sum_{i=1}^n \tau_i.
$$
These scheme appears in particular in the classical Monte - Carlo
method computing of definite integrals, may be multiple. \par

 \ In detail, let us consider for the problem of numerical computation
by the method Monte - Carlo the following definite (multiple, in
general case) integral
$$
 I := \int_D f(x) \ \mu(dx).
$$
Here $ \ \mu \ $ is a probability measure defined on the measurable
set  $ \ D: \ \mu(D) = 1. \ $  Let $ \  \{  \zeta_i \}, \ i =
1,2,\ldots,n, \ldots $ be a sequence of i., i.d. random variables
with distribution $\mu:$
$$
{\bf P}(\zeta_i \in A) = \mu(A)
$$
for all the measurable subsets $\{A\}$ of  the whole space $D$. \par

\vspace{3mm}

 \ The classical Monte - Carlo approximation $ \ \theta_n = I_n \ $ for
the source integral $I$ has the form
$$
 I_n := n^{-1} \sum_{i=1}^n f(\zeta_i), \ n \ge 2;
$$
i.e. here $\theta = I$. \par

On the other words, in this case

$$
  \tau_i = f(\zeta_i) - I.
$$
The consistent estimate $\hat{\beta}_n$ of the parameter $\beta$ as
$n \to \infty$, with the speed of convergence $\sim \ n^{-1/2} \
\ln^C n$ is offered for instance in the book \cite[chapter 5,
Section 5.12]{Ostrovsky1999}. \par

\vspace{4mm}

\noindent Define the (positive) value $ \ z^o(\delta), \ \delta \in
(0, 1/2) \ $ as follows
$$
 2^{2 - 2 z^o(\delta)} = \delta.
 $$
It follows immediately from the estimate of theorem 4.2a
 that with probability at least $1 - \delta$ and {\it for all
the values } $n \ge 2$
$$
|\theta_n - \theta| \le \frac{\beta \ \sqrt{2 \ \ln n }}{\sqrt{n}} \
\cdot ( 1 + z^o(\delta)).
$$

\vspace{4mm}

\noindent {\bf E.} Let us return to the mentioned above  article of Dominyka Kievinaite, Jonas Siaulys
\cite{Kievinaite - Saulys 2018}. One of the main result of one may be formulated as follows. Let $ \ \{\zeta_i\}, \ i = 1,2,3,\ldots, \ \zeta = \zeta_1 \ $ be a
sequence of i., i.d.  centered r.v.  belonging to certain  Grand Lebesgue Space; on the other words, satisfying  the famous Cramer's condition. \par
Let also $ \ d = \const > 0. \ $ Define the following probability

$$
V(x) = V_{\zeta}(x) =V_{\Law{\zeta}, d}(x) := {\bf P} \left( \ \sup_{n = 1,2,3,\ldots} \sum_{i=1}^n (\zeta_i - d) > x \right),\ x \ge 0.
$$

 \ It is proved in particular in \cite{Kievinaite - Saulys 2018} that under some additional conditions on the distribution $ \ \zeta \ $

$$
V(x) \le \min \left( \ 1,  c_1 \ e^{- c_2 x}  \   \right),\ c_1, c_2
=  c_1, c_2(\Law{\zeta}, d) = {\const} \in (0,\infty), \ x \ge 0.
$$

 \ In order to show that the last upper bound for this probability is essentially non - improvable,
 we bring a simple example. Assume that the r.v.  $ \ \zeta_0 \ $  is positive and has a standard
exponential distribution

$$
{\bf P}(\zeta_0 > x) = e^{-x}, \ x \ge 0.
$$
  \ Let us choose $ \ d = 1; \ $ notice that $ \  {\bf E} \zeta_0 = 1. \ $ Then the r.v.  $ \ \zeta_0 - 1 \ $ is centered and satisfies the Cramer's condition, as well as other suitable
conditions in the article  \cite{Kievinaite - Saulys 2018}. Wherein

$$
V_{\zeta_0}(x) \ge {\bf P}((\zeta_0 - 1) - 1 \ge x) = {\bf P}(\zeta_0 \ge 2 + x) = e^{-2} \cdot e^{-x}, \ x \ge 0.
$$

\vspace{4mm}

 \ {\bf Remark 8.1.  } \ It is interest to note that both the offered estimates are alike ones for obtained before {\it normalized}  sums of random variables. \par

\vspace{4mm}

\noindent {\bf F.} \ It is interest by our opinion  to obtain also bilateral  bounds for distribution of normed sums of weak or strong dependent random variables, for martingales etc. \par

\vspace{4mm}

\vspace{0.5cm} \emph{Acknowledgement.} {\footnotesize The first
author has been partially supported by the Gruppo Nazionale per
l'Analisi Matematica, la Probabilit\`a e le loro Applicazioni
(GNAMPA) of the Istituto Nazionale di Alta Matematica (INdAM) and by
Universit\`a degli Studi di Napoli Parthenope through the project
\lq\lq sostegno alla Ricerca individuale\rq\rq .\par

 \ The third author is grateful to M.Sgibnev for sending of Your interest article. \par


\begin{thebibliography}{100}

\bibitem{anatriellofiojmaa2015}
{\bf G.~Anatriello} and {\bf A.~Fiorenza.} {\it Fully measurable
grand Lebesgue spaces}. J. Math. Anal. Appl. \textbf{422} (2015),
no.~2, 783--797.

\bibitem{anatrielloformicaricmat2016}
{\bf G.~Anatriello} and {\bf M.~R.~Formica.} {\it Weighted fully
measurable grand Lebesgue spaces and the maximal theorem}. Ric. Mat.
\textbf{65} (2016), no.~1, 221--233.

\bibitem{Belyaev-Piterbarg1978} {\bf Yu.~K.~Belyaev} and {\bf V.~I.~Piterbarg.} {\it Random processes. Sample paths and intersections}. Collection of articles, Publishing House "MIR",
Moscow (1978); 249--257, (in Russian).


\bibitem{Bernstein1964}
{\bf S.N. Bernstein. } {\it  About some modifications of Chebyshev's inequality. }
Collected Works, (1964), Moscow, AN USSR, V. 4, pp. 330 - 337, (in Russian).

\bibitem{Braverman1991} {\bf M.~Sh.~Braverman.} {\it  Bounds on the sums of independent random variables in symmetric spaces.} Ukrainian Mathematical Journal,
 \textbf{43} (1991), no.~2, 148--153.

\bibitem{Braverman1994} {\bf M.~Sh.~Braverman.} {\it Independent Random Variables and Rearrangement Invariant Spaces.} London Mathematical Society, Lecture Notes Series, \textbf
{194}, Cambridge University Press, 1994.


\bibitem{Buldygin-Mushtary-Ostrovsky-Pushalsky} {\bf V.~V.~Buldygin, D.~I.~Mushtary, E.~I.
~Ostrovsky} and {\bf M.~I.~Pushalsky.} {\it New Trends in
Probability Theory and Statistics.} Mokslas (1992), V.1, 78--92;
Amsterdam, Utrecht, New York, Tokyo.

\bibitem{caponeformicagiovanonlanal2013}
{\bf C.~Capone, M.~R.~Formica} and {\bf R.~Giova.} {\it Grand
{L}ebesgue spaces with respect to measurable functions}. Nonlinear
Anal. \textbf{85} (2013), 125--131.



\bibitem{Ermakov etc. 1986}
{\bf S. V. Ermakov, and E. I. Ostrovsky.} {\it Continuity Conditions, Exponential Estimates, and the Central Limit Theorem for Random Fields.}
 Moscow, VINITY,  1986. (in Russian).



\bibitem{Fiorenza2000} {\bf A.~Fiorenza.} {\it Duality and reflexivity in grand Lebesgue
spaces.} Collect. Math. \textbf{51} (2000), no. 2, 131--148.

\bibitem{fiokarazanalanwen2004}
{\bf A.~Fiorenza} and {\bf G.~E.~Karadzhov.} {\it Grand and small
Lebesgue spaces and their analogs}, Z. Anal. Anwendungen \textbf{23}
(2004), no.~4, 657--681.

\bibitem{fioguptajainstudiamath2008}
{\bf A.~Fiorenza, B.~Gupta} and {\bf P.~Jain.} {\it The maximal
theorem for weighted grand Lebesgue spaces}. Studia Math.
\textbf{188} (2008), no.~2, 123--133.


\bibitem{Fiorenza-Formica-Gogatishvili-DEA2018}
{\bf A.~Fiorenza, M.~R.~Formica} and {\bf A.~Gogatishvili.} {\it On
grand and small Lebesgue and Sobolev spaces and some applications to
PDE's}. \emph{Differ. Equ. Appl.} \textbf{10} (2018), no.~1, 21--46.

\bibitem{fioforgogakoparakoNAtoappear}
{\bf A.~Fiorenza, M. R.~Formica, A.~Gogatishvili, T.~Kopaliani} and
{\bf J.~M. Rakotoson.} {\it Characterization of interpolation
between grand, small or classical Lebesgue spaces}. Preprint
arXiv:1709.05892, Nonlinear Anal., {to appear}.

\bibitem{fioformicarakodie2017}
{\bf A.~Fiorenza, M.~R.~Formica} and {\bf J.~M. Rakotoson.} {\it
Pointwise estimates for {$G\Gamma$}-functions and applications}.
Differential Integral Equations \textbf{30} (2017), no.~11-12,
809--824.

\bibitem{formicagiovamjom2015}
{\bf M.~R. Formica} and {\bf R.~Giova.} {\it Boyd indices in
generalized grand Lebesgue spaces and applications}. Mediterr. J.
Math. \textbf{12} (2015), no.~3, 987--995.

\bibitem{Hall-Heyde1980}
{\bf P.~Hall} and {\bf C.~C.~Heyde.} {\it Martingale limit theory
and its application. Probability and Mathematical Statistics}.
Academic Press, New York-London, 1980.

\bibitem{Hartman-Wintner1941}
{\bf P.~Hartman} and {\bf A.~Wintner.} {\it On the law of iterated
logarithm}. Amer. J. Math. \textbf{63}, (1941), 169--176.

\bibitem{Hu 1991}
{\bf T.~C.~Hu.} {\it  On the law of the iterated logarithm for
arrays of random variables}. Comm. Statist. Theory Methods
\textbf{20} (1991), no.~7, 1989--1994.

\bibitem{Hu-Weber 1992}
{\bf T.~C.~Hu} and {\bf N.~C.~Weber.} {\it On the rate of
convergence in the strong law of large numbers for arrays.} Bull.
Austral. Math. Soc., \textbf{45} (1992), no.~3, 479--482.

\bibitem{Iwaniec-Sbordone 1992}
{\bf T.~Iwaniec} and {\bf C.~Sbordone.} {\it On the integrability of
the Jacobian under minimal hypotheses.} Arch. Rational Mech. Anal.
\textbf{119} (1992), no.~2, 129--143.


\bibitem{Kievinaite - Saulys 2018}
{\bf  D. Kievinaite, J. Siaulys.} {\it \ Exponential bounds for the
tail probability of the supremum of an inhomogeneous random walk.}
Modern Stochastics: Theory and Applications, 5 (2), (2018), 129 \ -  \ 143. \\
\textbf{5.2} https://doi.org/10.15559/18-VMSTA99

\bibitem{Kolmogoroff 1929}
 {\bf A.~N.~Kolmogoroff.} {\it \"Uber das Gesetz des iterierten Logarithmus".} (German) Math. Ann., \textbf {101} (1929), no.~1, 126--135.

\bibitem{Kozachenko-Ostrovsky 1985}
{\bf Yu.~V.~Kozachenko} and {\bf E.~I.~Ostrovsky.} {\it The Banach
Spaces of random variables of sub-Gaussian type.} of Probab. and
Math. Stat., \textbf{32} (1985), (in Russian). Kiev, KSU, 43--57.

\bibitem{Kozachenko-Ostrovsky-Sirota Jan2017}
{\bf Yu.~V.~Kozachenko, E.~I.~Ostrovsky} and {\bf L.~Sirota.} {\it
Relations between exponential tails, moments and moment generating
functions for random variables and vectors.} arXiv:1701.01901v1
[math.FA], 8 Jan 2017.

\bibitem{Kozachenko-Ostrovsky-Sirota Oct2017}
{\bf Yu.~V.~Kozachenko, E.~I.~Ostrovsky} and {\bf L.~Sirota.} {\it
Equivalence between tails, Grand Lebesgue Spaces and Orlicz norms
for random variables without Cramer's condition.} arXiv:1710.05260v1
[math.PR],  15 Oct 2017.


\bibitem{Kozachenko at all 2018}
{\bf Yu.V. Kozachenko, Yu.Yu. Mlavets, and N.V. Yurchenko.} {\it Weak convergence of stochastic processes from spaces} $F_\psi(\Omega).$
STATISTICS, OPTIMIZATION AND INFORMATION COMPUTING, Vol.6,  June 2018, pp. 266 - 277.


\bibitem{liflyandostrovskysirotaturkish2010}
{\bf E.~Liflyand, E.~Ostrovsky} and {\bf L.~Sirota.} {\it Structural
properties of bilateral grand {L}ebesque spaces}. Turkish J. Math.
\textbf{34} (2010), no.~2, 207--219.


\bibitem{Naimark 2004}
{\bf Naimark B., Ostrovsky E.}
{\it Exact Constants in the Rosenthal Moment Inequalities for Sums of independent centered Random Variables.}
arXiv:math/0411614v1  [math.PR]  27 Nov 2004


\bibitem{Ng-Tang-Yang}
{\bf K.~W.~Ng, Q.~H.~Tang} and {\bf H.~Yang.} {\it Maxima of Sums of
Heavy-Tailed Random Variables.} ASTIN Bulletin: The Journal of the
IAA, {\textbf 32,} (2002), no.~1, pp. 43--55.


\bibitem{Ostrovsky1999}
{\bf E.~Ostrovsky.} {\it Exponential estimates for random fields and
its applications.} 1999, OINPE, Moscow - Obninsk.

\bibitem{Ostrovsky 1994}
{\bf E.~Ostrovsky.} {\it Exponential estimate in the Law of Iterated
Logarithm in Banach Space.} Math. Notes \textbf{56} (1994), no.~
5-6, 1165--1171.


\bibitem{Ostrovsky 2004}
{\bf E.~Ostrovsky.} {\it Bide-side exponential and  moment
inequalities for tails  of  distribution of  Polynomial
Martingales.} arXiv: math.PR/0406532 v.1 Jun. 2004.


\bibitem{Ostrovsky-Sirota Jan 2008}
{\bf E.~Ostrovsky} and {\bf L.~Sirota.} {\it  Exponential bounds in
the law of iterated logarithm for martingales.} arXiv:0801.2125v1
[math.PR], 14 Jan 2008.

\bibitem{Ostrovsky-Sirota-boundedeness operator bilateral GLS Oct2011}
{\bf E.~Ostrovsky} and {\bf L.~Sirota.} {\it Boundedness of
operators in bilateral Grand Lebesgue Spaces, with exact and weakly
exact constant calculation}, arXiv:1104.2963 [math.FA] Apr 2011.

\bibitem{Ostrovsky-Sirota Oct2013}
{\bf E.~Ostrovsky} and {\bf L.~Sirota} {\it Simplification of the
majorizing measures method, with development.} arXiv:1302.3202v1
[math.PR]  13 Feb 2013

\bibitem{Ostrovsky-Sirota Oct2015}
{\bf E.~Ostrovsky} and {\bf L.~Sirota.} {\it Vector rearrangement
invariant Banach spaces of random variables with exponential
decreasing tails of distributions.}
 arXiv:1510.04182v1 [math.PR] 14 Oct 2015.


\bibitem{Ostrovsky-Sirota-fundamental function GLS 2015}
{\bf E.~Ostrovsky} and {\bf L.~Sirota.} {\it Fundamental function
for Grand Lebesgue Spaces}, arXiv:1509.03644 [math.FA] Sept. 2015.


\bibitem{Pickands1967}
{\bf  J.~Pickands.} {\it  Maxima of stationary Gaussian processes.}
Z. Wahrscheinlichkeitstheorie und Verw. Gebiete \textbf{7} (1967),
190--223.

\bibitem{Pisier 1980}
{\bf G.~Pisier.} {\it Conditions d'entropie assurant la continuit\'e
de certains processus et applications \`{a} l'analyse harmonique.}
(French) Seminaire d'analyse fonctionnelle (1980), Exp. No. 13-14,
pp. 43 \ 46.

\bibitem{Qi 1994}
{\bf  Y.~C.~Qi.} {\it On strong convergence of arrays.} Bull.
Austral. Math. Soc. \textbf{50} (1994), no.~2, 219--223.

\bibitem{Rosalsky 1985}
{\bf  A.~ Rosalsky.} {\it On the number of successes in
independent trials.}  Sankhy\={a}, Ser. A \textbf{47} (1985),
no.~3, 380--391.

\bibitem{Samko-Umarkhadzhiev}
{\bf S.~G.~Samko} and {\bf S.~M.~Umarkhadzhiev.} {\it On
Iwaniec-Sbordone spaces on sets which may have infinite measure.}
Azerb. J. Math. \textbf{1} (1)\ (2011), 67--84.

\bibitem{Samko-Umarkhadzhiev-addendum}
{\bf S.~G.~Samko} and {\bf S.~M.~Umarkhadzhiev.} {\it On
Iwaniec-Sbordone spaces on sets which may have infinite measure:
addendum.} Azerb. J. Math \textbf{1} (2) \ (2011), 143--144.

\bibitem{Sgibnev 1996}
{\bf M.~S.~Sgibnev.} {\it On the distribution of the maxima of
partial sums.} Statist. Probab. Lett. \textbf{28} (1996), no.~3,
235--238.

\bibitem{Stout 1974}
{\bf W.~F.~Stout.} {\it Almost sure convergence.} Probability and
Mathematical Statistics, Vol. \textbf{24}. Academic Press. New
York-London, 1974.

\bibitem{Sung 1996}
{\bf  S.~H.~Sung.} {\it  An analogue of Kolmogorov's Law of the
Iterated  Logarithm for arrays.} Bull. Austral. Math. Soc.
\textbf{54} (1996),  no. 2, 177--182.

\bibitem{Teicher 1981}
{\bf  H.~Teicher.} {\it Almost certain behavior of row sums of
double arrays.} Analytical methods in probability theory
(Oberwolfach, 1980),  pp. 155--165, Lecture Notes in Math.,
\textbf{861}, Springer, Berlin-New York, 1981.

\bibitem{Welsch 1973}
{\bf R.~E.~Welsch.} {\it A convergence theorem for extreme values
from Gaussian sequences.} Ann. Probability  \textbf{1} (1973),
398--404.

\bibitem{Wittmann 1985}
{\bf 37. R.~ Wittmann.}  {\it  A General Law of Iterated
Logarithm.} Z. Wahrsch. Verw. Gebiete \textbf{68} (1985), no.~4,
521--543.







\end{thebibliography}
\end{document}